\documentclass[11pt]{amsart}
\usepackage{amsmath,amssymb,amsthm,amscd,latexsym}
\usepackage[all]{xy}
\usepackage[dvips]{graphicx}
\usepackage{hyperref}
\input xypic
\xyoption{all}

\newcommand{\rar}{\rightarrow}
\newcommand{\lar}{\longrightarrow}
\newcommand{\surjects}{\twoheadrightarrow}
\newcommand{\injects}{\hookrightarrow}

\newtheorem{Theorem}{Theorem}[section]
\newtheorem{Lemma}[Theorem]{Lemma}

\newtheorem{Proposition}[Theorem]{Proposition}
\newtheorem{Remark}[Theorem]{Remark}
\newtheorem{Example}[Theorem]{Example}

\newtheorem{Definition}[Theorem]{Definition}

\def\sqr#1#2{{\vcenter{\hrule height.#2pt
        \hbox{\vrule width.#2pt height#1pt \kern#1pt
            \vrule width.#2pt}
        \hrule height.#2pt}}}
\def\phi{\varphi}
\def\demo{\noindent{\bf Proof. }}
\def\square{\mathchoice\sqr64\sqr64\sqr{4}3\sqr{3}3}
\def\qed{\hspace*{\fill} $\square$}

\def\xx{{\bf x}}
\def\yy{{\bf y}}

\def\TT{{\bf T}}
\def\tt{{\bf t}}

\def\XX{{\mathbf X}}
\def\YY{{\bf Y}}

\def\ff{{\bf f}}

\def\ff{{\bf f}}

\def\hh{{\bf h}}

\def\fm{{\mathfrak m}}
\def\fn{{\mathfrak n}}
\def\NN{\mathbb N}

\def\hht{{\rm ht}\,}

\def\ass{{\rm Ass}\,}

\def\ker{{\rm ker}\,}

\def\rk{\rm rank}

\def\Rees{{\mathcal R}}
\def\cl#1{{\mathcal #1}}

\def\dd{{\mathbb D}}

\def\pp{{\mathbb P}}

\begin{document}

\title[An analogue of the Aluffi algebra]{An analogue of the Aluffi algebra for modules}%\\[3pt]

\author{Zaqueu Ramos}

\author{Aron Simis}

\keywords{Rees algebra, Aluffi algebra, symmetric algebra, module of derivations, normal module, Cohen--Macaulay, Gorenstein \\
\indent The first author was partially
supported by a CNPq post-doc fellowship (151229/2014-7).
\indent The second author was supported by a CNPq grant  (302298/2014-2) and a PVNS Fellowship from CAPES (5742201241/2016). He thanks the Departamento de Matem\'atica of the Universidade Federal da Paraiba for providing an appropriate environment for discussions on this work.\\
\indent Ramos address: Departamento de Matem\'atica, Universidade Federal de Sergipe, 49100-000 S\~ao Cristov\~ao, Sergipe, Brazil,
 email: zaqueu.ramos@gmail.com
\\
\indent Simis address: Departamento de Matem\'atica, Universidade Federal de Pernambuco, 50740-640 Recife, PE, Brazil, email: aron@dmat.ufpe.br
}

\bigskip

\begin{abstract}

P. Aluffi introduced in \cite{aluffi} a new graded algebra in order
to conveniently express characteristic cycles in the theory of
singular varieties. This algebra is attached to a surjective
ring homomorphism $A\surjects B$ by taking a suitable inverse limit
of graded algebras, one for each representation of $A$ as a residue
ring of a given ``ambient'' ring $R$. Since giving a ring surjection $A\surjects B$ is tantamount to giving an ideal $I\subset A$, it would seem natural to ask for an analogous notion for $A$-modules.
This is the central purpose of this work.
Since a given module may not admit any embedding into a free module, a preparatory toil includes dealing with this technical point at the outset.
On the bright side, the intrusion of modules raises a few algebraic questions interesting on their own.
It is to expect that this extension to modules may be transcribed in terms of coherent sheaves, thus possibly providing an answer to a question by Aluffi in this regard. Two main bodies of examples are treated in detail to illustrate how the theory works and to show the relation to finer properties of other algebras.

\end{abstract}

\maketitle

%\tableofcontents

\section*{Introduction}

In a seminal work (\cite{aluffi}), Aluffi introduced a so-called {\em quasi-symmetric} algebra as a direct limit of a family of algebras constructed from data depending on a variable ambient ring $R$.
One important property proved in \cite{aluffi} is that if for a given member of the family the ring $R$ is a regular local ring then the direct limit is identified with this particular member.
This motivated the study of one single such member in \cite{AluffiAA} and the corresponding algebra was dubbed an {\em embedded Aluffi algebra} -- we refer to the latter woek to avoid tedious repetition in this introduction.

From the last paragraph of \cite{aluffi} one reads:  ``Regarding a possible definition of quasi-symmetric algebras for coherent sheaves, this would presumably pivot on a good notion of Rees algebra of a module''.
The main goal of this paper is to pursue this suggestion, by drawing upon the existing wealth of material about Rees algebras of modules.
Our main reference for the latter is \cite{ram1} and the bibliography thereof.
We hope that the interested community of geometers moves in to transcribe the present algebraic account in terms of coherent sheaves to find suitable applications originally imagined in \cite{aluffi}.

%\smallskip

We now proceed to a short detailing of the sections.

The first section contains a brief review of the main facts about the Rees algebra of a module.
One of the subtleties of this notion is that there are various possible definitions, since an abstract module is not generally embeddable in a free module.
Thus, one has to cope with this virtual flaw of the theory and somehow move on to the case where all known definitions coincide.
For all purposes, one assumes throughout that the base ring is Noetherian.
Fortunately, as in the case of an ideal, when the module has a rank the dimension of its Rees algebra has a simple expression yielding much facilitation in the development.

The second section deals with the notion of the Aluffi algebra for a module. For its  precise formulation, the definition hinges on quite a bit of appropriate ingredients from module theory.
Albeit a certain instability in choosing these ground participants, the definition in the embedded case turns out to be an analog of the original one for ideals given in \cite{aluffi}.
To make sure this analog is a workable notion, one has to make sure of at least two things. First, that the embedded version admits a simple expression in terms of the Rees algebra of an embedded pull-back of the module.
Second and more difficult, that the embedded algebra lives in a family that allows to take inverse limits. This will in turn provide the analog of Aluffi's quasi-symmetric algebra for a module.

Most of the section is about meeting these two conditions.

The next two sections are dedicated to verifying the resulting impact of the theory so far on the case of two central modules in geometry: the module of vector fields (derivations) and the normal module of an ideal.

Thus, the third section deals with the module of derivations of a residue ring of a ground polynomial ring over a field. A natural embedded pull-back module is the so-called module of tangent vector fields.
This situation has been much studied in the case where the residue ring is the  coordinate ring of a  hypersurface.
In this setup, the module of tangent vector fields is also called {\em differential idealizer} (see \cite{s1}) in the algebraic category and {\em logarithmic derivations} (see \cite{Saito}) in the complex analytic category.

As a warm-up we first digress on the case where the hypersurface is a homogeneous free divisor. This is an ``easy'' case since the module of logarithmic derivations is a free module. As a consequence, the embedded Aluffi algebra of the module of derivations coincides with the corresponding symmetric algebra.
For convenience we give the detailed structure of an algebra presentation of the latter (Proposition~\ref{free_divisor}).

The main result of the section is Theorem~\ref{smooth_hypersurface}.
It completely describes the homological nature of the (embedded) Aluffi algebra in the case of a smooth projective hypersurface.
This result is based on two essential other results, Theorem~\ref{Rees_of_Z} and Theorem~\ref{Pfaff-Gorenstein}.
The first of these is a vast generalization of a previous result proved independently by Herzog--Tang--Zarzuela (\cite[Corollary 2.3]{HTZ}) and Goto--Hayasaka--Kurano--Nakamura (\cite[Theorem 4.1]{GHKN}).
The second is a recent result of A. Kustin (\cite{Kustin}),  responding affirmatively to a question of ours.
An alternative path for our query had been proposed earlier to us by S. Goto (\cite{Goto}), framed on a general view of Gorenstein ideals of codimension one in Cohen--Macaulay rings.

The main homological feature in this result totally crumbles down if the hypersurface has singularities, even in the plane case.

The last section discusses the embedded Aluffi algebra of the normal module to an ideal in a polynomial ring over a field.
This algebra has no interest in the case where the ideal is a complete intersection since the normal module is free.
The section deals with the next simplest case, namely, that of a codimension $2$ perfect ideal $J$ which is presented by linear syzygies.
Here the main result is Theorem~\ref{Main_Normal}.

Alas, in the case of the normal module, the Aluffi algebra will fail to be a Gorenstein, or even Cohen--Macaulay, ring.
Yet, to our surprise, the case of a codimension $2$ linearly presented perfect ideal $J$ turns out to really be a result about linearly presented modules of projective dimension one, with some additional features. A subset of the results in the theorem draws on the treatment in \cite{ram1}.
In addition, we have also used the fact that, in the generic case, the ideal $J$ is rigid and that it is a {\em licci} ideal.
Under the latter hypothesis, a crucial specialization result takes place, recently communicated to us by B. Ulrich.

By and large, the more structured parts of Theorem~\ref{Main_Normal} are controlled by subtle conditions on the local number of generators,  and by the relation between the dimension of the ambient polynomial ring and the number of generators of the ideal $J$.

\medskip

{\sc Acknowledgements.} We are greatly indebted to S. Goto, A. Kustin and B. Ulrich for their assistance in crucial parts of the work; the second author also thanks C. Miranda Neto for some exchange in an early exposition of the topics.

\section{Generalities on symmetric and Rees algebras of modules}

In this section we review for the reader's convenience some definitions and general
 facts about Rees algebras of modules, as taken from \cite{ram1} (cf. also \cite{elam99}).

Let $A$ be a Noetherian ring and let $E$ be an $A$--module.
Throughout $E$ will be finitely generated, but many module-theoretic notions do not require this proviso.

A torsion element of $E$ is an element $z\in E$ annihilated by some
nonzerodivisor of $A$. The set of such elements forms a submodule
of $E$, denoted $\tau_A(E)$, and called the $A$-{\it torsion
submodule\/} of $E$, or simply the $A$-torsion of $E$. Clearly,
$\tau_A(E)=\ker (E\rar E\otimes _AK)$, where $K$ is the total ring
of quotients of $A$. $E$ is said to be torsionfree if
$\tau_A(E)=0$. Clearly, $E/\tau_A(E)$ is torsionfree.

We define the {\it Rees algebra\/} $\Rees (E)$ of $E$ to be the
residue algebra of the symmetric algebra ${\mathcal S}(E)$ by its
$A$-torsion\index{torsion} $\cl T$. Note that ${\mathcal S}(E)$ is
naturally an $\NN$-graded $A$-algebra and that $\cl T$ is actually
a homogeneous ideal in this grading. The module $E$ is said to be
of {\it linear type\/} if  ${\mathcal T}=0$, i.e., if the canonical
map ${\mathcal S}(E)\surjects \Rees( E)$ is an isomorphism.

Given a free presentation of a finitely generated $A$-module $E$
$$G\lar F \lar E \rar 0,$$
the symmetric algebra has a  presentation ${\mathcal S}(E)\simeq {\mathcal S}(F)/{\rm Im}(G){\mathcal S}(F)$.
In a sufficiently concrete translation, it is the residual $A$-algebra of a polynomial ring $A[\mathbf{Y}]$ over $A$ by the ideal $I_1(\mathbf{Y})\cdot Z)$, where $Z$ is a representative matrix of the map $G\rar F$ -- in other words, the so-called {\em syzygies} of $E$ in the given presentation.
Clearly, then the Rees algebra of $E$ has a presentation $\Rees (E)\simeq A[\mathbf{Y}]/(I_1(\mathbf{Y})\cdot Z), \widetilde{\cl T})$, where $\widetilde{\cl T}$ is the corresponding torsion lifted to the polynomial ring.

 If $(A,\fm)$ is Noetherian local with maximal ideal $\fm$ or standard graded over a field with irrelevant ideal $\fm$, the module $E$ (or its Rees algebra) is said to be of {\em fiber type} if $\widetilde{\cl T}$ is generated by forms in $(A/\fm)[\mathbf{Y}]$.

A finitely generated $A$-module $E$ is said to {\it
have a {\rm (}generic {\rm )} rank}, say, $r$, if $E_p$ is $A_p$-free of rank $r$ for
every prime $p\in \ass A$. If this is the case, we denote $\rk
E=r$. Equivalently, if $K$ is the total ring of quotients of $A$,
then $E$ has rank $r$ if and only if $K \otimes_A E \cong K^r$.

One says that $E$ is {\it torsionless\/} if the natural
$A$-homomorphism $E\rar E^{**}$ into the double $A$-dual is
injective. Equivalently, $E$ is torsionless if and only it can be
embedded into a free $A$-module. Clearly, a torsionless module is
torsionfree (by the latter characterization). Conversely, it is
not difficult to see that if $E$ is a finitely generated torsionfree $A$-module having a rank then it is torsionless.

We remark that if $E$ is a submodule of a free $A$-module $G$,
many authors define the Rees algebra of $E$ to be the image of the
natural map ${\mathcal S}(E) \rar {\mathcal S}(G)$ (see \cite{Rees} for
such an approach).  The two definitions coincide when the finitely
generated module $E$ is torsionfree and has a rank (i.e., the
$A$-map $E\rar K\otimes _A E$ is injective and $K\otimes _A E$ is
$K$-free) as in this situation the kernel of ${\mathcal S}(E) \rar
{\mathcal S}(G)$ is the $A$-torsion submodule of ${\mathcal S}(E)$. The
present definition depends solely on $E$, not on any particular
embedding of $E$. Finally, in particular, note that if $I\subset
A$ is an ideal containing a regular element (i.e., $I$ has rank
one) then the present definition retrieves the ordinary notion,
namely, the graded $A$-algebra $A\oplus I\oplus
I^2\oplus\cdots\simeq A[It]\subset A[t]$.

One ought to mention more encompassing approaches to the notion of the Rees algebra of a module. The priority is \cite[Chapitre III]{Mic}, where a fairly satisfactory theory is established in the case of a finitely generated module over a Noetherian ring.
A more refined notion is given in \cite{EHU}. Both theories eventually show that the definition boils down to the case of a submodule of a free module. There is a relevant overlapping between the two ways of reasoning and it might turn out to be useful to explain this in more detail.
The definition in both approaches, as well the one used in this paper, all coincide in the case of a finitely generated module with a rank.

\medskip

In the case where the given module has a rank, one has a pretty  dimension formula for $\Rees (E)$.
\begin{Proposition} \label{2.2} {\rm (\cite{ram1})}
Let $A$ be a Noetherian ring of
dimension $d$ and $E$  a finitely generated $A$--module having a rank
$r$. Then $ \dim \Rees(E) = d+r = d + {\rm height} ~~ \Rees(E)_+.$
\end{Proposition}

The Rees algebra $\Rees (E)$ has a kind of universal mapping property parallel to that
of ${\mathcal S}(E)$, to wit,
given a torsionfree $A$-algebra $B$ and an $R$-module homomorphism $E\rar B$ there exists a unique
$A$-algebra map $\Rees (E)\rar B$ making the obvious diagram commutative.
Yet some care has to be exercised when dealing with torsion as it behaves quite badly with respect to ring change $A \rar A'$. As a result, one of the
drawbacks of the Rees algebra, in contrast to the symmetric
algebra, is that it is not truly functorial in the sense of commuting with ring  change.

If $(A,\fm)$ is Noetherian local with maximal ideal $\fm$ or standard graded over a field with irrelevant ideal $\fm$, the $A/\fm$-algebra $\Rees(E)/\fm\Rees(E)$ is called the {\em special fiber} of $\Rees(E)$ and its dimension is the {\em analytic spread} of $E$.

\smallskip

While the symmetric algebra is a more na\"ive algebraic notion, not so much a formula for its dimension.
The fundamental theoretic result is the formula of Huneke--Rossi, showing that the dimension of the symmetric algebra of a finitely generated module $E$ over a Noetherian ring $A$ is given by the so-called Forster number
$$\sup _{p\in {\rm Spec}(A)}\,\{\dim A/p+\mu_p(E)\},$$
where $\mu_p(E)$ denotes the minimal number of generators of $E_p$ as an $A_p$-module.

If $(A,\fm)$ is local or standard graded over a field ($E$ graded) then we get $\dim S(E)\geq \mu(E)$ (global minimal number of generators).

The above formulas or bounds are not effective.
For an effective version of formulas in the case the module has a rank, we refer to \cite{SV1}, \cite{SV2}.
There one introduces the numerical invariant
$${\rm df}(E)=\max_{1\leq t\leq \rk(\phi)}\{\max\{0, \rk(\phi)-t+1- {\rm cod}(I_t(\phi))\},$$
associated to a free presentation of $E$ with syzygy matrix $\phi$.

Then, provided $A$ is equidimensional and catenarian, one has
$$\dim S(E)= \dim A+\rk(E)+ {\rm df}(E).$$

In addition, in the case $E$ has a rank, one says that $E$ satisfies the condition $(F_0)$ if the Fitting deffect ${\rm df}(E)$ vanishes, in which case one often says that the symmetric algebra $S(E)$ has the expected dimension (often an unfortunate terminology).
A stronger assumption, with no impact on the dimension level, is the condition $(F_1)$ (also called $G_{\infty}$ in the Artin--Nagata terminology). The latter has an impact on the question of the irreducibility of ${\rm Proj}(S(E))$.

The following sort of weakening is useful:  a finitely generated $A$-module $E$ having rank $r$ is said to satisfy condition $G_s$, for a given integer $s\geq 1$, if $\mu (E_{\mathfrak p})\leq \dim A_{\mathfrak p}+r-1$ whenever $1\leq \dim A_{\mathfrak p}\leq s-1$.
Note that, if $A$ is local or standard graded over a field, such a condition for $s=\dim A$ affects only the behavior of $E$ on the punctured spectrum of $A$.

We refer to \cite{BHbook} for the remaining current algebraic notions to be used in this paper, mainly those concerning graded structures, such as the properties of being Cohen--Macaulay, Gorenstein and so forth.

\section{The construction}

Let $A$ be a Noetherian ring and let $E$ be a finitely generated
$A$-module. The following notion of rank avoids the
restriction that $E$ have a generic rank.
\begin{Definition}\rm The {\it embedding rank\/} of $E$ is the least integer $s\geq 1$
for which there is an embedding $E/\tau_A(E)\hookrightarrow F$,
where $F$ is free of rank $s$. If $E/\tau_A(E)$ is not torsionless then no such embedding exists, in which case one may say
that $E$ has no embedding rank (or its embedding rank is undefined).
\end{Definition}
As observed, $E$ admits an embedding rank if and only if the
torsionfree module $E/\tau_A(E)$ is torsionless. It is also clear
that the embedding rank is just the rank of $E$ if it happens to
have one. As an example, a nonzero ideal $I\subset A$ without a
rank (i.e., contained in some associated prime $P$ of $A$)
still has embedding rank one.

Now, let $R\surjects A$ be a surjective ring homomorphism and let
$J=\ker (R\surjects A)$. If ${\bf F}$ is an $R$-free module of
rank $n$ and $F$ is an $A$-free module of rank $n$, an $R$-module
surjection ${\bf F}\surjects F$ is said to be {\it compatible\/}
with the structural ring surjection $R\surjects A$ if $\ker ({\bf
F}\surjects F)=J{\bf F}$.

Let $E$ be a finitely generated $A$-module such that $E/\tau_A(E)$
is torsionless (but may not have a rank). Fix an embedding
$E/\tau_A(E)\hookrightarrow F$, where $F$ is $A$-free of rank $n$.
\begin{Definition}\rm
A $J$-{\it deformation\/} of (the given embedding of) $E$ to $R$
is a finitely generated $R$-module ${\bf E}$ such that:
\begin{enumerate}
\item[{\rm (i)}] There is an embedding ${\bf E}/\tau_R({\bf E})\hookrightarrow {\bf F}$
where ${\bf F}$ has rank $n$ and, moreover, $J\subset ({\bf
E}/\tau_R({\bf E}))\colon {\bf F}$
\item[{\rm (ii)}] For any surjective $R$-map ${\bf F}\surjects F$ compatible with the
ring map $R\surjects
A$, the restriction to ${\bf E}/\tau_R({\bf E})$ maps onto
$E/\tau_A(E)\subset F$.
\end{enumerate}
\end{Definition}

We will say that two $J$-deformations ${\bf E}$ and ${\bf E}_1$ of
the given embedding of $E$ to $R$ are {\it congruent\/} if there
is an $R$-isomorphism ${\bf F}\simeq {\bf F}_1$ of the respective
ambient free $R$-modules which maps ${\bf E/\tau_R({\bf E})}$ onto
${\bf E}_1/\tau_R({\bf E}_1)$ and which commutes with the
respective structural $R$-module surjections ${\bf F}\surjects F$
and  ${\bf F}_1\surjects F$ (i.e., which maps $J{\bf F}$ to $J{\bf
F}_1$).

Note that the conditions imply in particular that, for a given
$J$-deformation ${\bf E}$ to $R$, there is an $A$-module
isomorphism $({\bf E}/\tau_R({\bf E}))/J{\bf F}\simeq E/\tau_A(E)$
and that if ${\bf E}_1$ is another $J$-deformation
 to $R$ then ${\bf E}_1/\tau_R({\bf E}_1)\simeq {\bf
E}/\tau_R({\bf E})$.

As to the existence, fix a free $R$-module ${\bf F}$ of rank $n$
and fix an $R$-map surjection ${\bf F}\surjects F$ compatible with
$R\surjects A$. Consider the inverse image of $E/\tau(E)\subset F$
by ${\bf F}\surjects F$. This is clearly a $J$-deformation of $E$
to $R$. Therefore, any finitely generated module $E$ with  a fixed
embedding $E/\tau(E)\hookrightarrow F$ admits a unique
$J$-deformation ${\bf E}$ to $R$ up to $R$-torsion and congruence.

Note the induced $R$-algebra surjection ${\mathcal S}_R({\bf
E}/\tau({\bf E}))\surjects {\mathcal S}_A(E/\tau_A(E))$ of symmetric
algebras.

We can now state the analogue of Aluffi's notion.

\begin{Definition}\rm Let $E$ be a finitely generated $A$-module such that $E/\tau_A(E)$
is torsionless. Fix an embedding $E/\tau_A(E)\hookrightarrow F$,
where $F$ is $A$-free of rank $n$. Let $R\surjects A$ be a ring
surjection with kernel $J\subset R$. Let ${\bf E}$ denote a
torsionfree $J$-deformation of $E$ to $R$ (which, by the
preceding, is uniquely defined up to congruence). The {\it {\rm
(}$R$-embedded{\rm)} Aluffi algebra\/} of $E$ is
$${\mathcal A}_{_{R\surjects A}}(E)\colon ={\mathcal S}_A(E/\tau_A(E))\otimes_{{\mathcal S}_R({\bf E})}
\Rees_R({\bf E}).$$
\end{Definition}

\begin{Remark}\rm This treatment extends
Aluffi's original definition to finitely generated modules and
besides does not assume a priori that $E$ is torsion free, merely
that $E/\tau_A(E)$ has a positive embedding rank $r$  and a fixed
embedding into a free $A$-module of rank $r$ is given. Of course, if
$E$ is actually a nonzero ideal ${\mathfrak a}\subset A$ (i.e., an
$A$-module with embedding rank one with a fixed embedding) then the
present notion retrieves Aluffi's concept.
\end{Remark}

\begin{Lemma}\label{aluffi2rees}
There is an $A$-algebra isomorphism
\begin{equation}\label{presentation}
{\mathcal A}_{_{R\surjects A}}(E)\simeq \frac{\Rees_R({\bf E})}{(J,J{\bf F})\,\Rees_R({\bf
E})},
\end{equation}
where $J$ is in degree $0$ and $J{\bf F}$ is in degree
$1$. In particular, there is a surjective $A$-algebra homomorphism
${\mathcal A}_{_{R\surjects A}}(E)\surjects \Rees_A(E)$.
\end{Lemma}
\demo  By the very definition of ${\mathcal A}_{_{R\surjects A}}(E)$,
it is equivalent to show that ${\mathcal S}_A(E/\tau_A(E))\simeq {\mathcal
S}_R({\bf E})/(J,J{\bf F})\,{\mathcal S}_R({\bf E})$ (recall that
$\tau_R({\bf E})={0}$ by construction). But this follows from the
ordinary universal property of the symmetric algebra.

To see the surjection, one considers any algebra in sight in terms of a
presentation over $R$. For this, we lift a set of generators of
$E$ over $A$ to vectors in ${\bf E}$. These lifted vectors
together with the set of natural generators of the submodule
$J{\bf F}$ induced from a set of generators of $J$ will be a set
of generators of ${\bf E}$. We accordingly take a set of
presentation variables $\TT=\{\TT\setminus\TT'\}\cup \TT'$, where
the variables $\TT'$ correspond to the generators coming from
$J{\bf F}$. Let then $\Rees_R({\bf E})\simeq R[\TT]/{\mathcal J}$ and
$\Rees_A(E)\simeq R[\TT\setminus\TT']/(J,\widetilde{{\mathcal J}})$,
where $\widetilde{{\mathcal J}}$ denotes the lifting to
$R[\TT\setminus\TT']$ of a set of generators of the presentation
ideal of $\Rees_A(E)$ over $A[\TT\setminus\TT']$. By the above
isomorphism it is clear that ${\mathcal A}_{_{R\surjects A}}(E)\simeq
R[\TT]/(J,\TT',{\mathcal J})\simeq
R[\TT\setminus\TT']/(J,\overline{{\mathcal J}})$, where
$\overline{{\mathcal J}}$  is generated by the generators of ${\mathcal
J}$ by setting $\TT'=0$. Therefore, one has an inclusion of ideals
$(J,\overline{{\mathcal J}})\subset (J,\widetilde{{\mathcal J}})$, as
required. \qed

\medskip

As in \cite{aluffi}, one considers the inverse system of  ring surjections $R\rar R'$
commuting with given fixed surjections $R\surjects A$ and
$R'\surjects A$. Let $E$ be a finitely generated $A$-module such
that $E/\tau_A(E)$ is torsionless, with a fixed embedding
$E/\tau_A(E)\hookrightarrow F$, where $F$ is $A$-free of rank $n$.
Consider the (congruence classes of ) torsionfree $J$-deformations
${\bf E}$, one for each $J=\ker (R\surjects A)$. Clearly, any
surjection $R\rar R'$ in the inverse system induces a
natural map ${\bf F}\rar {\bf F}'$ that maps ${\bf E}$ to ${\bf
E}'$ and is compatible with the structural surjections ${\bf
F}\surjects F$ and  ${\bf F}'\surjects F$. The resulting set is
again an inverse system inducing an
inverse system of $A$-algebra surjections ${\mathcal A}_{_{R\surjects
A}}(E)\rar {\mathcal A}_{_{R'\surjects A}}(E)$. The resulting inverse limit
$$\cl A _A(E)\colon = \lim_{_{R\surjects A}}\,{\mathcal A}_{_{R\surjects A}}(E)$$
is the (absolute) {\it Aluffi algebra\/} of a finitely generated $A$-module $E$ with
a fixed embedding into a
free $A$-module.

\bigskip

Still as in \cite{aluffi}, we next restrict our inverse system to the category of
algebras of finite type over a ground field $k$. The main bulk of the theory is
given by the following theorem.

\begin{Theorem}\label{reduction_to_regular} Let $A$ be an algebra of finite type over a perfect field $k$ and let $E$ be a
finitely generated $A$-module with a fixed embedding into a free $A$-module. Let
$R\surjects A$ be a $k$-algebra surjection, where $R$ is a smooth algebra of finite
type over $k$. Then the structural map of $A$-algebras $\cl A _A(E)\surjects {\mathcal
A}_{_{R\surjects A}}(E)$ is an isomorphism.
\end{Theorem}
\demo
Although the proof follows the one of \cite[Theorem 2.9]{aluffi}, by
retro-tracing the argument to a situation of rings and ideals, the details of the argument are more delicate. We first slightly rephrase the ring setup in a somewhat more precise form. Thus, let $R$ and $R'$ be smooth algebras of
finite type over $k$ with given surjective $k$-algebra homomorphism $\pi: R\surjects
R'$ commuting with fixed surjective $k$-algebra homomorphisms $R\surjects A$ and
$R'\surjects A$. According to the remarks prior to the statement there is a surjective $A$-algebra
homomorphism ${\mathcal A}_{_{R\surjects A}}(E)\stackrel{\tilde{\pi}}{\surjects} {\mathcal
A}_{_{R'\surjects A}}(E)$ and as in the proof of  \cite[Theorem 2.9]{aluffi} it
suffices to show that this map is an isomorphism.

To do this, one uses the auxiliary commutative diagram of $k$-algebra homomorphisms
$$\begin{array}{ccc}
&R\otimes_k R'&\\
\kern10pt\raise5pt\hbox{$\iota$}\nearrow && \kern-8pt\searrow \raise5pt\hbox{$\mu$}\\
\kern-5ptR & \stackrel{\pi}{\lar\kern-12pt\lar} & \kern5pt R'\\
\kern10pt\searrow &  &\kern-8pt \swarrow\\
&A&
\end{array}
$$
where $\iota$ is the one-sided injective map $a\mapsto a\otimes 1$ and $\mu$ is the
multiplication map $a\otimes a'\mapsto \pi(a)a'$, while the lower maps are the
structural maps. Note first that $\mu$ is split by the other one-sided injection
$R'\simeq 1\otimes_k R'\subset R\otimes_k R'$. Also, one has $R\otimes_k
R'/(\dd,1\otimes J')\simeq A$, where $J'=\ker(R'\surjects A)$ and $\dd$ is the
kernel of $\mu$.

One now claims that the diagram induces a diagram of $A$-algebra homomorphisms
$$\begin{array}{ccc}
&{\mathcal A}_{_{R\otimes_k R'\surjects A}}(E)&\\
\kern12pt\raise5pt\hbox{$\tilde{\iota}$}\nearrow && \kern-20pt\searrow
\raise5pt\hbox{$\tilde{\mu}$}\\
\kern-5pt {\mathcal A}_{_{R\surjects A}}(E) & \stackrel{\tilde{\pi}}{\lar\kern-12pt\lar} &
\kern5pt {\mathcal A}_{_{R'\surjects A}}(E)\\
\end{array}
$$
where $R\otimes_k R'\surjects A$ is the composite map $R\otimes_k
R'\stackrel{\mu}{\rar} R'\surjects A$ and $\tilde{\mu}$ is a split $A$-algebra
homomorphism.

We first deal with this latter map. By definition,  given an embedding
$E/\tau_A(E)\subset F$, we choose a map ${\bf F}_{R'}\surjects F$ compatible with
the structural surjection $R'\surjects A$ and take the inverse image ${\bf
E}_{R'}\subset {\bf F}_{R'}$, so that ${\bf E}_{R'}/J'{\bf E}_{R'}\simeq
E/\tau_A(E)$. Next, we apply to ${\bf F}_{R'}$ base change by the one-sided
inclusion $R'\injects R\otimes_kR'$ and set ${\bf F}_{R\otimes_kR'}=R\otimes_k {\bf
F}_{R'}$. Thus, we get a commutative diagram where the horizontal $A$-module maps are split maps
$$\begin{array}{ccccc}
{\bf E}_{R'} & \injects & R\otimes_k {\bf E}_{R'} & \surjects &  {\bf
E}_{R'}\\
\bigcap &&\bigcap && \bigcap\\
{\bf F}_{R'} & \injects & R\otimes_k {\bf F}_{R'} & \surjects & R\otimes_k {\bf
F}_{R'}
\end{array}
$$
It follows that $R\otimes_k {\bf E}_{R'}$ is a $(\dd,1\otimes J')$-deformation of $E/\tau_A(E)$. Thus, set ${\bf E}_{R\otimes_kR'}:=R\otimes_k {\bf E}_{R'}$. Next,
since $R\otimes_k {\mathcal S}_{_{R'}}({\bf E}_{R'}) \simeq {\mathcal S}_{_{R\otimes_kR'}}
({\bf E}_{R\otimes_kR'})$ and $R\otimes_k {\mathcal R}_{_{R'}}({\bf E}_{R'})\simeq {\mathcal
R}_{_{R\otimes_kR'}} ({\bf E}_{R\otimes_kR'})$, this diagram induces a commutative
diagram of symmetric and Rees algebras

$$\begin{array}{cccccc}
{\mathcal R}_{_{R'}}({\bf E}_{R'}) & \injects & {\mathcal R}_{_{R\otimes_kR'}}
({\bf E}_{R\otimes_kR'}) & \stackrel{\rho}{\surjects} &  {\mathcal R}_{_{R'}}({\bf E}_{R'})&\\
\uparrow &&\uparrow  && \uparrow &\\[-10pt]
\uparrow &&\uparrow  && \uparrow &\\
{\mathcal S}_{_{R'}}({\bf E}_{R'}) & \injects & {\mathcal S}_{_{R\otimes_kR'}} ({\bf
E}_{R\otimes_kR'}) & \surjects &  {\mathcal S}_{_{R'}}({\bf E}_{R'}) &\surjects {\mathcal
S}_{_{A}}(E/\tau_A(E))
\end{array}
$$
where the only non-trivial map is $\rho$. The existence of this map is guaranteed by \cite[Proposition 4.5]{elam99} because $R\otimes_k {\bf E}_{R'}$ is locally of linear type at the associated primes of ${\mathbb D}=\ker(R\otimes_kR'\surjects R')$ since $R'$ is locally regular (being smooth over a perfect field).

One then gets a split injection of Aluffi algebras

$${\mathcal S}_{_{A}}(E/\tau_A(E))
\otimes_{_{ {\mathcal S}_{_{R'}}({\bf E}_{R'})}} \left({\mathcal R}_{_{R'}}({\bf
E}_{R'})\right) \injects {\mathcal S}_{_{A}}(E/\tau_A(E))\otimes_{_{ {\mathcal
S}_{_{R\otimes_kR'}} ({\bf E}_{R\otimes_kR'})}} \left({\mathcal R}_{_{R\otimes_kR'}}
({\bf E}_{R\otimes_kR'})\right)
$$
whose left inverse is precisely the required map $\tilde{\mu}$.

On the other hand, one has

\begin{eqnarray*}&&{\mathcal S}_{_{A}}(E/\tau_A(E))\otimes_{_{ {\mathcal S}_{_{R\otimes_kR'}}
({\bf E}_{R\otimes_kR'})}} \left({\mathcal R}_{_{R\otimes_kR'}}
({\bf E}_{R\otimes_kR'})\right)\\
&\simeq & {\mathcal S}_{_{A}}(E/\tau_A(E))\otimes_{_{R\otimes_k {\mathcal S}_{_{R'}}({\bf
E}_{R'})}} \left(R\otimes_k {\mathcal R}_{_{R'}}({\bf E}_{R'})\right)\\
&\simeq & \left({\mathcal S}_{_{A}}(E/\tau_A(E))\otimes_k R\right) \otimes_{_{ {\mathcal
S}_{_{R'}}({\bf E}_{R'})}}
\left({\mathcal R}_{_{R'}}({\bf E}_{R'})\right)\\
&\simeq & \left({\mathcal S}_{_{A}}(E/\tau_A(E))\otimes_k A\right) \otimes_{_{ {\mathcal
S}_{_{R'}}({\bf E}_{R'})}}
\left({\mathcal R}_{_{R'}}({\bf E}_{R'})\right)\\
&\simeq & {\mathcal S}_{_{A}}(E/\tau_A(E)) \otimes_{_{ {\mathcal S}_{_{R'}}({\bf E}_{R'})}}
\left({\mathcal R}_{_{R'}}({\bf E}_{R'})\right)
\end{eqnarray*}
It follows that $\tilde{\mu}$ is an isomorphism.

It remains to explain the map $\tilde{\iota}$ and argue that it is injective.
This will conclude the proof.

As before, we construct everything {\it ab initio}. This time around, we choose a
map ${\bf F}_{R}\surjects F$ compatible with the structural surjection $R\surjects
A=R/J$ and take the inverse image ${\bf E}_{R}\subset {\bf F}_{R}$, so that ${\bf
E}_{R}/J{\bf E}_{R}\simeq E/\tau_A(E)$. Next, we apply to ${\bf F}_{R}$ base change
by the one-sided inclusion $R\injects R\otimes_kR'$ and set ${\bf F}_{R\otimes_kR'}=
{\bf F}_{R}\otimes_k R'$. Note this induces a surjection ${\bf
F}_{R\otimes_kR'}\surjects {\bf F}_{R}\otimes_R R'/J'={\bf F}_{R}\otimes_R A=F$
compatible with the structural map $R\otimes_kR'\surjects A$.

We note that though $R\injects R\otimes_kR'$ has no left inverse, it is an injective
flat map because $k$ is a field and, for any $R$-module $M$ one has
$M\otimes_R(R\otimes_kR')\simeq M\otimes_k R'$ as $R\otimes_kR'$-modules. We claim
that ${\bf E}_{R}\otimes_k R'$ is a $(\dd,1\otimes J')$-deformation of
$E/\tau_A(E)$. Since $J{\bf E}_{R}\subset {\bf F}_{R}$, it follows that
$(\dd,1\otimes J'){\bf E}_{R}\otimes_k R'\subset {\bf F}_{R}\otimes_k R'$. This shows
that condition (i) of the definition of deformation is satisfied.

To see the second condition of this definition, let $\psi:{\bf F}_{R}\otimes_k
R'\surjects F$ be a map compatible with $R\otimes_kR'\surjects
A=R\otimes_kR'/(\dd,1\otimes J')$. Since the image of ${\bf E}_{R}$ in ${\bf
E}_{R}\otimes_k R'$ generates the latter as $R\otimes_kR'$-module and the
restriction of ${\bf F}_{R}\surjects F$ to ${\bf E}_{R}$ maps onto $E/\tau_A(E)$ -
by condition (ii) applied to the $J$-deformation ${\bf E}_{R}$ - it is also clear
that the restriction of $\psi:{\bf F}_{R}\otimes_k R'\surjects F$ to ${\bf
E}_{R}\otimes_k R'$ maps onto $E/\tau_A(E)$ as well, so condition (ii) also holds
here.

Having got so far, the conclusion is easy. Indeed, first we have
$${\mathcal S}_{_{R\otimes_kR'}}({\bf E}_{R}\otimes_k R')\simeq {\mathcal S}_{_{R}}
({\bf E}_{R})\otimes_{_{R}}(R\otimes_kR') \simeq {\mathcal S}_{_{R}}({\bf
E}_{R})\otimes_{_{k}}R'
$$
and similar isomorphisms for the corresponding Rees algebra. Therefore, one has
\begin{eqnarray*}
&&{\mathcal S}_{_{A}}(E/\tau_A(E))\otimes_{_{ {\mathcal S}_{_{R\otimes_kR'}} ({\bf
E}_{R\otimes_kR'})}} \left({\mathcal R}_{_{R\otimes_kR'}}
({\bf E}_{R\otimes_kR'})\right)\\
&\simeq& {\mathcal S}_{_{A}}(E/\tau_A(E))\otimes_{_{ {\mathcal S}_{_{R}}({\bf E}_{R})
\otimes_{_{k}}R'}} \left({\mathcal R}_{_{R}}({\bf E}_{R})
\otimes_{_{k}}R'\right)\\
&\simeq& \left({\mathcal S}_{_{A}}(E/\tau_A(E))\otimes_{_{ {\mathcal S}_{_{R}}({\bf E}_{R})
}} {\mathcal R}_{_{R}}({\bf E}_{R}) \right)\otimes_{_{k}} R'
\end{eqnarray*}
Since the map
$${\mathcal S}_{_{A}}(E/\tau_A(E))\otimes_{_{ {\mathcal S}_{_{R}}({\bf E}_{R})
}} {\mathcal R}_{_{R}}({\bf E}_{R})\lar \left({\mathcal S}_{_{A}}(E/\tau_A(E)) \otimes_{_{
{\mathcal S}_{_{R}}({\bf E}_{R}) }} {\mathcal R}_{_{R}}({\bf E}_{R}) \right)\otimes_{_{k}}
R'
$$
is injective, we are through.
\qed

A major point has to do with the dimension of the Aluffi algebra.
For simplicity, we assume that $A$ is local and that $E$ has a rank.
\begin{Proposition}\label{dimensionbounds}
Let $R$ be a smooth finitely generated algebra over a perfect field $k$ and let  $A$ be a catenarian equidimensional ring with
$A\simeq R/J$,  $\dim A=d$ and $\hht (J)=g$.
Let $E$ be a finitely generated $A$-module with rank $r$ and a
given embedding $E/\tau_A(E)\subset F$ into a free $A$-module of
rank $n$. Then
\begin{equation}\label{dimensioninequalities}
\min\{d+n+g,\, d+r+{\rm df}(E)\}\geq \dim \cl A _{A}(E)\geq d+r.
\end{equation}
\end{Proposition}
\demo The only mysterious number here is ${\rm df}(E)\}$, which has been explained in the preliminaries of Section 1.
Thus, the left inequality follows immediately from the definitions since by the
latter the Aluffi algebra $\cl A _A(E)$ is a residue algebra both of ${\mathcal
S}_A(E/\tau_A(E))$ and of $\Rees_R({\bf E})$.
The right inequality follows at once from the surjective
$A$-algebra homomorphism in Lemma\,\ref{aluffi2rees}. \qed

\section{The Aluffi algebras of two classical modules}

The previous sections made it possible to study the Aluffi algebra as a single member of its defining direct limit, as long as the ambient $R$ is a finitely generated smooth algebra over a perfect field.

In this and later sections we assume that, in fact,
 $R=k[x_1,\ldots,x_n]$ is a polynomial ring over a field $k$ of characteristic zero
and let $J\subset R$ be an ideal.
We will be concerned with two basic modules in both commutative algebra and algebraic geometry: the module of $k$-derivations (vector fields) of $R/J$ and the normal module of $J$ in $R$.
An advantage at the outset is that both modules are torsion free $R/J$-modules by definition, being $R/J$-duals of the well-known module of K\"ahler differentials and the conormal module of $J$ in $R$, respectively.

Alas, the complete understanding of the Aluffi algebra of any of these modules in general is a task far from immediate grasp.
We will study this algebra for particular cases of the ideal $J$.

\subsection{The module of vector fields of a hypersurface}

Set $E:={\rm Der}_k(R/J)$, the module of $k$-derivations of the $k$-algebra $R/J$.
Note that $E\subset F:=\sum_{i=1}^n (R/J)\,\partial/\partial x_i$, so in particular it is a torsionfree
$R/J$-module.
Then $E$ has a deformation to the {\em module of logarithmic derivations} $$\mathbf{E}={\rm Der}_J(R):=\{ \delta\in {\rm Der}_k(R)\;|\; \delta (J)\subset J\}.$$
This module is a  global polynomial analogue of Saito's module ${\rm Der}_S({\rm log}\,V(J))$ in case $S=\mathbb{A}^n$ and $V(J)$ is a hypersurface (\cite[(1.4) Definition]{Saito}).
Note that the module of logarithmic derivations sits naturally as an $R$-submodule of $\mathbf{F}:=
\sum_{i=1}^n R\,\partial/\partial x_i={\rm Der}_k(R)$.

In this environment the presentation (\ref{presentation}) of the Aluffi algebra is in general still out of reach due to
the mysterious structure of the Rees algebra of the module of logarithmic derivations.
However, assuming $J$ reduced, one has the dimension bounds
\begin{equation}\label{dim_idealizer_case}
2n-2\,\hht(J)\leq\dim \cl A _{R/J}({\rm Der}_k(R/J))\leq \dim R +
\rk ({\rm Der}_J(R))-1= 2n-1
\end{equation}
from the general results in the previous sections.

\smallskip

We will focus on the case where $J$ is a principal ideal generated by a reduced polynomial $f$.
In this situation, by (\ref{dim_idealizer_case}) the dimension of the corresponding Aluffi algebra gets trapped between $2n-2$ and $2n-1$.

On a pinhead, the two main ingredients to ascribe a presentation of the Aluffi algebra is a presentation of the Rees algebra $\mathcal{R}_R(\mathbf{E})$ of the deformation $\mathbf{E}$ and a set of generators of $f\mathbf{F}$ viewed as a subideal.

Now, $\mathcal{R}_R(\mathbf{E})$ can be reached in two ways since $\mathbf{E}\subset \mathbf{F}$. The first is as the algebra obtained from the symmetric algebra of $\mathbf{E}$ by killing $R$-torsion.

The second is as follows: trading notation $\partial/\partial x_i=t_i$ for visual convenience, the symmetric algebra of $\mathbf{F}$ is identified
with the polynomial ring $R[t_1,\ldots, t_n]$ and the Rees algebra of the module of logarithmic derivations of $f$ gets identified with the $R$-subalgebra
generated by the $1$-forms which are the entries of the matrix product $\mathbf{t}\cdot \mathcal{D}$, where $\mathbf{t}=(t_1\cdots t_n)$
and $\mathcal{D}$ denotes the matrix whose columns generate the module of logarithmic derivations, i.e., we get
$$\Rees ({\rm Der}_f(R))\simeq R[\mathbf{t}\cdot \mathcal{D}]\subset R[\mathbf{t}].$$

In the second approach to further identify $f\mathbf{F}$ as  contained in the graded component of degree one
$R[\mathbf{t}\cdot \mathcal{D}]_1$, one  needs to know how $f$ sits, in the original module notation, as an element in the $i$th row of $ \mathcal{D}$, for some $i=1,\ldots,n$.

\subsubsection{The homogeneous case}

Let $R=k[\xx]=k[x_1,\ldots,x_n]$ denote a standard graded polynomial ring over a field $k$.
Let $f\in R=k[\xx]=k[x_1,\ldots,x_n]$ stand for a reduced form of degree $d\geq 2$ and let $\mathbf{E}:={\rm Der}_f(R)\subset
{\rm Der}(\mathbf{F})=\oplus_i R\,\partial/\partial x_i\simeq R^n$ and $E:={\rm Der}_k(R/(f))\subset
\oplus_i (R/(f))\,\partial/\partial x_i\simeq (R/(f))^n$
denote the $R$-module of logarithmic derivations and the module of derivations of $R/(f)$, respectively.
Note that both modules are graded submodules of the respective graded free modules in the grading inherited from
the standard grading of $R$.
Since both are generated in varying degrees, the respective associated $R$-algebras (symmetric, Aluffi, Rees, special
fiber) are quasi-homogeneous (i.e., weighted homogeneous). In particular, some care has to be exercised when using
their Hilbert series in order to deduce ideal theoretic properties.

\subsubsection*{Homogeneous free divisors}

For a free divisor $f\in R$ the Aluffi
algebra of the module of derivations of $R/(f)$ coincides with its symmetric algebra. We include it anyway for the sake of completion.

Recall that $f$ is said to be a free divisor (\`a la Saito) if the module of logarithmic derivations ${\rm Der}_f(R)$ is a free module (necessarily, with a basis of $n$ vectors).

\begin{Proposition}\label{free_divisor}
Let $f\in R$ be a homogeneous free divisor. Then
\begin{enumerate}
\item[{\rm (a)}] $\cl A _{R/(f)}(E)\simeq \cl S _{R/(f)} (E)$  % $\simeq R[U_1,\ldots, U_n]$ {\rm (}polynomial ring{\rm )}
\item[{\rm (b)}] Let $\cl A _{R/(f)}(E)\simeq R[U_1,\ldots, U_n]/\mathfrak{A}$ denote a presentation
of the Aluffi algebra of $E$ as a quasi-homogeneous $R$-algebra, where $U_1,\ldots,U_n$ are weighted indeterminates over $R$ with
$R[U_1,\ldots, U_n]\simeq \cl R _R(\mathbf{E})$ by mapping $U_i$ to a corresponding basis element of $\mathbf{E}$
in an ordered basis.
Then $\mathfrak{A}$ is a codimension $2$ perfect ideal whose presentation matrix is the $(n+1)\times n$ matrix
$$\left(
  \begin{array}{ccc}
    U_1 & \cdots & U_n \\
    \hline \\
     & \boldsymbol{\cl E} &  \\
  \end{array}
\right),
$$
where $\boldsymbol{\cl E}$ denotes the matrix whose columns are an ordered basis of $\mathbf{E}$ {\rm (}up to signs{\rm )}.
\item[{\rm (c)}] The presentation matrix of $E={\rm Der}_k(R/(f))$ as an $R$-module of projective dimension $1$ is the $n\times n$
matrix of co-factors of $\boldsymbol{\cl E}$.
\end{enumerate}

\end{Proposition}
\demo
(a) This is trivial since $\mathbf{E}$ is free.

(c) This matrix is the content matrix of $f\mathbf{F}$ in $\mathbf{E}$.
Since $f=\det \boldsymbol{\cl E}$ (up to a sign perhaps), the assertion is clear.

(b) This is the core.
As control steps, note that dimensions are right.
Thus, with $\mathfrak{A}$ of codimension $2$, the virtual dimension is $\dim \cl A _{R/(f)}(E)=2n-2$.
This is the true dimension if $(F_0)$ holds over $R/(f)$.
Now the rank of $E$ is $\dim R/(f)=n-1$; by tensoring in (c) the presentation matrix of $E$ over $R/(f)$ is
still the matrix of co-factors modulo $(f)$ and its rank over $R/(f)$ is $1$,
On the other hand, the cofactors
generate an ideal of codimension $2$ over $R$, hence of codimension $1$ over $R/(f)$.
Thus, $(F_0)$ holds trivially ($1=\rk -t+1=1-1+1=1$).

The proof is also easy.
Indeed, by the structural form in Lemma~\ref{aluffi2rees}, $\mathfrak{A}$ is generated
by $f$ and the entries of the various product matrices $(U_1...U_n)\cdot C_i$, where $C_i$ stands for the content
matrix of $(0\,\ldots\, 0\, f\, 0\, \ldots\, 0)^t$ ($f$ in the $i$th position) in $\boldsymbol{\cl E}$.
By (c) we know that the entries of $C_i$ are suitable co-factors.
\qed

\subsubsection*{Smooth projective hypersurfaces}

The following preliminary proposition substantially extends a result of  Herzog--Tang--Zarzuela (\cite[Corollary 2.3]{HTZ}) and Goto--Hayasaka--Kurano--Nakamura (\cite[Theorem 4.1]{GHKN}).

\begin{Theorem}\label{Rees_of_Z}
Let $\ff=\{f_1,\ldots,f_n\}$ be a regular sequence of forms of same degree $d\geq 1$ in $R=k[\xx]=k[x_1,\ldots,x_n], (n\geq 2)$.
Consider the cokernel
\begin{equation}\label{the_Z}
\bigwedge^3 R^n \stackrel{\kappa}{\lar} \bigwedge^2 R^n \lar \mathbf{Z}\rar 0,
\end{equation}
where $\kappa=\kappa(\bf f)$ denotes the second map in the Koszul complex of $\ff$.
Then the Rees algebra $\cl R_R(\mathbf{Z})$ of the $R-$module $\mathbf{Z}$ is a Gorenstein domain of codimension ${{n-1}\choose {2}}$ on the polynomial ring $R[y_{ij}\,|\,{1\leq i<j\leq n}]$.
\end{Theorem}
\demo
Let $\YY:=(y_{ij})_{1\leq i<j\leq n}$ denote a generic $n\times n$ skew-symmetric matrix over $R$.
Consider the augmented skew-symmetric matrix
\begin{equation}\label{matrix}
\mathcal{P}=
\begin{pmatrix}
\YY & ^t \ff \\
-\ff & 0
\end{pmatrix},
\end{equation}
where we have written $\ff=(f_1\,\cdots \,f_n)$.

We will show that the defining ideal of  $\cl R_R(\mathbf{Z})$ on the polynomial ring $R[y_{ij}\,|\,{1\leq i<j\leq n}]$ is the ideal $\mathfrak{P}$ generated by the $4$-Pfaffians of (\ref{matrix}) and that this ideal has the required properties in the statement.

\smallskip

In the completely generic case -- i.e., if the augmented column has indeterminate entries mutually independent from the $y_{ij}$'s -- it is fairly classical that the ideal of the $4$-Pfaffians is Gorenstein, of codimension
${{n-1}\choose {2}}$, generated by ${{n+1}\choose {4}}$ quadrics
\cite{Kleppe}.
In the present case we have the same number of generators, ${{n}\choose {4}}$ of the $4$-Pfaffians are quadrics, while the remaining ${{n}\choose {3}}$ ones are the $4$-Pfaffians involving last column and row, hence are of degree $d+1$.
Explicitly, these equations have the following respective shapes:

\begin{equation}
\label{inner_pfaffians}
y_{ik}y_{jl} - y_{jk}y_{il} - y_{ij} y_{kl} , 1 \leq i < j <k <l \leq n
\end{equation}
and
\begin{equation}
\label{outer_pfaffians}
f_iy_{jk} - f_j y_{ik} + f_ky_{ij} , 1 \leq i < j <k \leq n.
\end{equation}

\smallskip

We now prove that the defining ideal of  $\cl R_R(\mathbf{Z})$ is the ideal $\mathfrak{P}$ generated by the $4$-Pfaffians of the matrix (\ref{matrix}).
We may assume that $d\geq 2$ as otherwise one has  $n$ independent linear forms, so by a change of variables this is the case of Remark~\ref{previous}.

\medskip

{\sc Claim 1:} The ideal $\mathfrak{P}$ is Gorenstein of codimension ${{n-1}\choose {2}}$.

{\sc Proof:} Let $\XX=\{X_1,\ldots, X_n\}$ be algebraically independent over $R[\yy]$.
Consider the completely generic version of (\ref{matrix}):
\begin{equation} \label{matrix_generic}
%\widetilde{\mathcal{P}}=
\begin{pmatrix}
\YY & ^t \XX \\
-\XX & 0
\end{pmatrix}.
\end{equation}
Clearly, $\XX-\ff=\{X_i-f_i, 1\leq i\leq n\}$ is a regular sequence on $R[\XX]=k[\xx,\XX]$.

In addition, $\XX-\ff$ is also a regular sequence on $R[\YY,\XX]/\widetilde{\mathfrak{P}}$, where  $\widetilde{\mathfrak{P}}$ is the ideal of $4$-Pfaffians of (\ref{matrix_generic}).
To see this one can use a Gr\"obner basis sort of argument as follows.
Fix the graded lex order in the standard graded polynomial ring $k[\xx,\YY,\XX]$ with the following ordering of variables: $\xx>\YY>\XX$.
Recall that we are assuming that  $\ff$ is generated in degree $d\geq 2$ in $k[\xx]$; then there is an inclusion of initial ideals ${\rm in} (\ff)\subset {\rm in}(\XX - \ff)$.

Indeed, first it is clear that the initial term of any $f_i$ is the initial term of the corresponding $X_i-f_i$.
Next, consider the $S$-polynomial $g_{i,j}$ of a pair $f_i, f_j$, for different indices $i,j$.
Let the respective leading terms be denoted $\fn_i$ and $\fn_j$.
We may as usual assume that $\gcd(\fn_i,\fn_j)\neq 1$; say, $e:=\deg(\gcd(\fn_i,\fn_j))$.
Consider in parallel the $S$-polynomial $h_{i,j}$ of the pair $X_i-f_i, X_j-f_j$.
Then $h_{i,j}=g_{i,j}+ T$, where $T$ is certain combination of $X_i,X_j$ with coefficients in $k[\xx]$.
Note that any nonzero term of  $g_{i,j}$  has degree $d+d-e=2d-e$. At the other end, any nonzero term of $T$ has degree $d-(e-1)$.
Since $2d-e>d-(e-1)$ for $d\geq 2$ and the order is graded  then the leading term of $g_{i,j}$ and of  $h_{i,j}$ are the same.

Next consider the $S$-polynomial of  $g_{i,j}$ and some $f_l$ (possibly, $l=i$ or $l=j$).
By the same token, the nonzero terms of this polynomial have degree higher that the degrees of any nonzero term in the part involving any variable form $k[\XX]$.
This argument goes on till one eventually reaches a Gr\"obner basis of $\ff$, thus proving our intermediate contention.

Now we apply this as follows:

$$\left( {\rm in}(\widetilde{\mathfrak{P}}), {\rm in}(\ff)\right)\subset \left( {\rm in}(\widetilde{\mathfrak{P}}),{\rm in}(\XX-\ff)\right)\subset  {\rm in} \left(\widetilde{\mathfrak{P}}, \XX-\ff \right)
$$

Since the leftmost ideal is a sum of two ideals in separate variables, its codimension is the sum of the respective codimensions, i.e., equals ${{n-1}\choose {2}}+n$.
Therefore, the codimension of the rightmost ideal is bounded below by this number.
This proves that  $\XX-\ff$ is also a regular sequence on $R[\YY,\XX]/\widetilde{\mathfrak{P}}$.
Therefore, specializing by $\XX-\ff$ yields that the specialized ideal  $\mathfrak{P}$ is Gorenstein of same  codimension as that of $\widetilde{\mathfrak{P}}$.

\begin{Remark}\rm
The above seems  to work under more general assumptions.
We suspect that it suffices to assume that the ideal $(\ff)$ has codimension at least $n-2$. However, the above proof breaks down at the end.
\end{Remark}

%We argue as in the proof of \cite[Proposition 3.1]{Kleppe}.

\medskip

{\sc Claim 2:} The ideal $\mathfrak{P}$ is prime.

\smallskip

{\sc Proof:} Since $\mathfrak{P}$ is homogeneous it suffices to prove that it is a normal ideal.
Since it is Gorenstein (hence, Cohen--Macaulay), it certainly satisfies Serre's property $(S_2)$.
Therefore, it remains to prove that it satisfies property $(R_1)$.

Set $N:={n-1\choose 2}$, the codimension of $\mathfrak{P}$.
Let $\Theta$ denote the Jacobian matrix of the generators of $\mathfrak{P}$ with respect to all variables in sight.
The goal is to show that
\begin{equation}\label{Serre_R1}
{\rm cod}(I_{N}(\Theta)+I)\geq N+2,
\end{equation}
where $I_r(M)$ denotes the ideal generated by the $r$-minors of a matrix $M$.
A close inspection of the generators of $\mathfrak{P}$ readily shows that
$$\Theta=\left(\begin{array}{c|c}{\bf 0}&A\\\hline B&C\end{array}\right)$$
for certain submatrices $A,B$ and $C$, where the latter is of order ${n\choose 3}\times {n\choose 2}$  and each one of its entries is either some $f_i$ ($1\leq i\leq n$) or else zero.

In addition, one has the following property: $(f_1^N,\ldots,f_n^{N})\subset I_N(C).$
To see this, one needs to go into further detail about $C$ by noting that $C$ is the Jacobian matrix of the forms
$$h_{ijk}:=f_iy_{jk}-f_jy_{ik}+f_ky_{ik}
\quad(1\leq i<j<k\leq n),$$
with respect to the $\yy$ variables.
Let us observe that, given $1\leq \ell\leq n,$ there are exactly $N$ among these forms $h_{ijk}$ such that $\ell\in\{i,j,k\}$ -- let $\hh_{\ell}$ denote this set.
Looking this way, one sees that the Jacobian matrix of $\hh_{\ell}$ with respect to the variables $y_{ij}$ such that $\ell\notin \{i,j\}$
is a diagonal submatrix of $C$ with $\pm f_{\ell}$ throughout the main diagonal.
Clearly then $f_{\ell}^N\in I_N(C)$, showing the contention.

Applying this inclusion, one has
\begin{eqnarray} \nonumber
{\rm cod}(I_N(\Theta)+I)&\geq& {\rm cod}({\rm Pf}_4(\yy),f_1^N,\ldots,f_n^N)\\ \nonumber
&=&{\rm cod}({\rm Pf}_4(\yy))+{\rm cod}(f_1^N,\ldots,f_n^N)\\ \nonumber
&=&{n-2\choose 2} +n\\ \nonumber
&=& N+2,
\end{eqnarray}
as claimed in (\ref{Serre_R1}).
\qed

\medskip

To conclude the proof of item (ii), we first note the codimension of the Rees algebra:
\begin{eqnarray}\nonumber
{\rm cod}(\cl R_R(\mathbf{Z}))&=&\dim R[\yy]-\dim \cl R_R(\mathbf{Z})\\ \nonumber
&=&n+{{n}\choose {2}}- (\dim R + \rk\, \mathbf{Z})\\ \nonumber
&=&
{{n}\choose {2}}-(n-1)={{n-1}\choose {2}}
\end{eqnarray}

In order to close the argument it suffices as a consequence of Claim 1 and Claim 2 to show that $\mathfrak{P}$ is contained in the presentation ideal of $\cl R_R(\mathbf{Z})$.
For this, we use the following device in order to replace the generators of $\mathbf{Z}$ by their corresponding ``polynomializations''.
Namely, as an $R$-subalgbra of $R^n$, the Rees algebra $\cl R_R(\mathbf{Z})$ is also isomorphic to the $R$-subalgebra
of the polynomial ring $R[\tt]=R[t_1,\ldots,t_n]$ generated by the biforms
$$\{f_it_j-f_jt_i,  1 \leq i < j\leq n\}$$
coming from the Koszul relations of $\ff$.
An immediate inspection yields that the two sets (\ref{inner_pfaffians}) and (\ref{outer_pfaffians}) of equations vanish for $y_{ij}\mapsto f_it_j-f_jt_i$, with $ 1 \leq i < j\leq n$.

This completes the proof.
\qed

\begin{Remark}\label{previous}
\rm
The proof that the defining ideal of  $\cl R_R(\mathbf{Z})$ is the ideal $\mathfrak{P}$ generated by the $4$-Pfaffians of the matrix (\ref{matrix}) has been proved in \cite[Corollary 2.3]{HTZ} and \cite[Theorem 4.1]{GHKN} in the case where $\ff=\{x_1,\ldots,x_n\}$, assuming char$(k)\neq 2$.
In this situation, the argument has been based on detecting a Sagbi basis (in the case of \cite{HTZ}) or drawing upon generalized Grassmann algebras (in the case of \cite{GHKN}).
Examining our line of argument above, is becomes apparent that in this case all one needs is that the ideal of $4$-Pfaffians in the generic case be prime, a result due to Kleppe \cite{Kleppe}.
\end{Remark}

\begin{Theorem}\label{smooth_hypersurface}
Let $f\in R=k[\xx]=k[x_1,\ldots,x_n]$ be a form of degree $d\geq 2$ defining a smooth projective hypersurface, where $k$ has characteristic larger than $\deg(f)$.
Keep the remaining notation as above.
Then:
\begin{enumerate}
\item[{\rm (i)}] The Rees algebra $\cl R_R(\mathbf{E})$ is a polynomial ring in one variable $T$ over the Rees algebra $\cl R_R(\mathbf{Z})$ of the cokernel
\begin{equation}\label{The_Z_bis}
\bigwedge^3 R^n \stackrel{\kappa(\partial f)}{\lar} \bigwedge^2 R^n \lar \mathbf{Z}\rar 0,
\end{equation}
where $\kappa(\partial f)$ denotes the second map in the Koszul complex of the partial derivatives of $f$.
\item[{\rm (ii)}] The Rees algebra $\cl R_R(\mathbf{E})$ is a Gorenstein domain of codimension ${{n-1}\choose {2}}$ on the polynomial ring $R[y_{ij}\,|\,{1\leq i<j\leq n}][T]$.
\item[{\rm (iii)}] The Aluffi algebra $\cl A _{R/(f)}(E)$ is a Gorenstein ring of codimension
${{n-1}\choose {2}}+2$ on the polynomial ring $R[y_{ij}\,|\,{1\leq i<j\leq n}][T]${\rm;} in particular, its dimension coincides with the dimension of the Rees algebra $\cl R_{R/(f)}(E)$ of the module of derivations over $R/(f)$.
\end{enumerate}
\end{Theorem}
\demo
(i) The assumption on the characteristic of $k$ implies, via Euler's formula, that $f$ belongs to the ideal of the partial derivatives of $f$, in which case the module of logarithmic derivations $\mathbf{E}$ is a direct sum of the first module of syzygies $\mathbf{Z}$ of the gradient ideal of $f$ and the cyclic module $R\epsilon$ generated by the Euler derivation $\epsilon$ (\cite[Proposition 3.2]{s1}).
Since taking the symmetric algebra is functorial, one has
$$\cl S_R(\mathbf{E})\simeq \cl S_R(\mathbf{Z})\otimes_R \cl S_R (R\epsilon)\simeq
\cl S_R(\mathbf{Z})\otimes_R R[T]\simeq\cl S_R(\mathbf{Z})[T],$$
where $T$ is an indeterminate over $R$.
On the other hand, $\cl R_R(\mathbf{E})\simeq \cl S_R(\mathbf{E})/\tau_R(\cl S_R(\mathbf{E}))$, where $\tau_R(\cl S_R(\mathbf{E}))$ is the $R$-torsion of  $\cl S_R(\mathbf{E})$. Because taking torsion commutes with a free extension, we get
$$\cl R_R(\mathbf{E})\simeq \cl S_R(\mathbf{Z})[T]/\tau_R(\cl S_R(\mathbf{Z})[T])\simeq (\cl S_R(\mathbf{Z})/\tau_R(\cl S_R(\mathbf{Z}))[T]\simeq \cl R_R(\mathbf{Z})[T].$$
To conclude, since $V(f)$ is smooth, the partial derivatives of $f$ generate an $(x_1,\ldots,x_n)$-primary ideal, hence form a regular sequence.
This implies that the module $\mathbf{Z}$ of their syzygies coincides with the cokernel of the second Koszul map
$$\bigwedge^3 R^n \stackrel{\kappa(\partial f)}{\lar} \bigwedge^2 R^n.$$

\smallskip

(ii) By (i), this an immediate consequence of Theorem~\ref{Rees_of_Z}.

\smallskip

(iii) Set $\YY:=\{y_{ij},|, 1\leq i<j\leq n\}$.
By (i) and the proof of Theorem~\ref{Rees_of_Z} one has
$$\cl R_R(\mathbf{E})\simeq R[\YY,T]/\left(\mathfrak{P}\right),$$
where $\mathfrak{P}={\rm Pff}_4(\mathcal{P})$ and $\mathcal{P}$ is as in (\ref{matrix}).

Therefore, according to the presentation (\ref{presentation}) in Lemma~\ref{aluffi2rees}, one has
\begin{equation}
\label{presentation2}
\cl A _{R/(f)}(E)\simeq R[\mathbf{Y}, T]/\left(\mathfrak{P}, f, f\mathbf{F}\right),
\end{equation}
where $f\mathbf{F}\subset {\mathbf E}=\mathbf{Z}+R\epsilon$, with $\epsilon$ denoting the Euler derivation.

While the Rees algebra of $\mathbf{E}$ has a standard graded structure, the Aluffi algebra $\cl A _{R/(f)}(E)$ has a nonstandard graded structure over $R[\mathbf{Y}, T]$ if $\deg(f)\geq 3$, with $T$ of weight $1$ (for the Euler derivation $\epsilon$) and each $y_{ij}$ of weight $\deg(f)-1$.

With these weights, the defining equations coming from the inclusion
$$f\mathbf{F}\subset {\mathbf E}=\mathbf{Z}+R\epsilon$$
become homogeneous.
These equations are written in compact form as
\begin{equation}\label{key_ideal}
I_1\left((\mathbf{Y}, T)\cdot \mathcal{C}\right),
\end{equation}
with $\mathcal{C}$ the transpose of the concatenation $[\kappa\,|\,^t\partial f]$, where $\kappa$ is the usual matrix of the first Koszul map of the variables $x_1,\ldots,x_n$  and
$$\partial f:=\left(\frac{\partial f}{\partial x_1}\, \cdots \, \frac{\partial f}{\partial x_n}\right).$$
Collecting (\ref{presentation2}) and (\ref{key_ideal}) yields
\begin{equation}
\label{presentation3}
\cl A _{R/(f)}(E)\simeq R[\mathbf{Y}, T]/\left(\mathfrak{P}, f, I_1\left((\mathbf{Y}, T)\cdot \mathcal{C}\right)\right),
\end{equation}

{\sc Claim.} $\left(f, I_1\left((\mathbf{Y}, T)\cdot \mathcal{C}\right)\right)=\left(I_1\left((\xx, T)\cdot \mathcal{P}\right)\right)$

\smallskip

This is an expression of the so-called Jacobian duality (\cite[Section 1]{dual}) plus the easy observation that $f$ is essentially the product of $(\xx, T)$ by the rightmost column of $\mathcal{P}$.

\smallskip

As a consequence, one has the presentation $\cl A _{R/(f)}(E)\simeq R[\mathbf{Y}, T]/(\mathfrak{P},I_1(\xx, T)\cdot \mathcal{P})$
 where the presentation ideal is generated by the defining ideal of the symmetric algebra of the cokernel of $\mathcal{P}$ and the $4$-Pfaffians of $\mathcal{P}$.

 \smallskip

 We now consider the generic analogue of this ideal, namely, replace $\mathcal{P}$ by the generic skew-symmetric matrix $\mathcal{G}$ of the same size -- thus, the partial derivatives of $f$ are to be replaced by new variables $y_{1,n+1},\ldots, y_{n,n+1}$.

 \smallskip

 In this regard, the following result has just been proved by A. Kustin:

 \begin{Theorem}{\rm (\cite[Theorem 4.8]{Kustin})}\label{Pfaff-Gorenstein}
The ideal $(\mathfrak{G},I_1(\xx, T)\cdot \mathcal{G}))$ is Gorenstein of codimension ${{n-1}\choose {2}}+2$ on the polynomial ring $R[y_{ij}\,|\,{1\leq i<j\leq n+1}][T]$, where $\mathfrak{G}={\rm Pff}_4(\mathcal{G})$.
 \end{Theorem}

We will specialize from the above generic result. For this it suffices to prove that the  sequence $\{y_{1,n+1}-\partial f/\partial x_1,\ldots, y_{n,n+1}-\partial f/\partial x_n\}$ is regular on
the polynomial ring $R[y_{ij}\,|\,{1\leq i<j\leq n+1}][T]$ and also modulo the ideal $(\mathfrak{G},I_1(\xx, T)\cdot \mathcal{G}))$.

The first assertion is clear. To show the second assertion it is enough to argue that the ideal $(\mathfrak{P},I_1(\xx, T)\cdot \mathcal{P})$ has codimension ${n-1\choose 2}+2$ on $R[y_{ij}\,|\,{1\leq i<j\leq n+1}][T]$.
To see this, proceed as follows.

First, by (\ref{dim_idealizer_case}) one has
\begin{equation*}
2n-2\leq\dim \cl A _{R/(f)}(E)\leq  2n-1.
\end{equation*}
Therefore, it suffices to show that ${\rm cod}\, (\mathfrak{P},I_1(\xx, T)\cdot \mathcal{P})\geq {n-1\choose 2}+2.$
At this point we follow a similar argument as in the proof of \cite[Lemma]{Kustin}, namely, let $P\subset R[y_{ij}\,|\,{1\leq i<j\leq n+1}][T]$ be a prime ideal containing $(\mathfrak{P},I_1(\xx, T)\cdot \mathcal{P}).$

If $T\in P$ then $Q:=(\mathfrak{P},T)$ is a prime of codimension ${n-1\choose 2}+1$ contained in $P.$ But, since $f\in P\setminus Q$ then  ${n-1\choose 2}+2\leq \hht P.$

If $T\not\in P$,  let $\mathcal{P}'$ and $\mathcal{P}''$ denote the matrices obtained from $\mathcal{P}$ by removing, respectively, its last column and its last row and last column.
Set $\mathfrak{P}'':={\rm Pf}_4(\mathcal{P}'').$
Observe  the following facts:  $\mathfrak{P}''$ is a prime ideal of codimension ${n-2\choose 2}$ contained in $P$ and the entries of $(\xx,T)\cdot \mathcal{P}'$ form a regular sequence on  $R_T/\mathfrak{P}''R_T$ contained in $PR_T.$

It follows that $${n-1\choose 2}+2={n-2\choose 2}+n\leq {\rm cod} \,PR_T={\rm cod}\, P,$$
as was to be shown.
\qed

\begin{Remark}
\rm
In addition to the result of the last theorem, one has that the $R/(f)$-module $E={\rm Der}_k(R/(f))$ is locally free on the punctured spectrum of $R/(f)$;
in particular, the defining ideal of the Rees algebra $\cl R_{R/(f)}(E)$ is the saturation by $(\xx)$ of the defining ideal of the Aluffi algebra $\cl A _{R/(f)}(E)$.
To see this one argues as follows: since $V(f)$ is smooth, the module of Ka\"hler $k$-differentials
of $R/(f)$ is locally free on the punctured spectrum of $R/(f)$ and hence so is its $R$-dual $E={\rm Der}_k(R/(f))$.
Therefore, $\cl A _{R/(f)}(E)\simeq \cl R_{R/(f)}(E)$ locally on the punctured spectrum of $R/(f)$, which proves the supplementary statement on the saturation.
\end{Remark}

 This is as much as one can go for the comparison of the two algebras: simple examples show that the Rees algebra of ${\rm Der}_k(R/(f))$, although Cohen--Macaulay, is not Gorenstein.
From the other end, for an arbitrary singular $f$ there does not seem to be a clear pattern.
Thus, e.g., for the cuspidal cubic $f$ in $\pp^2$ the Rees algebra $\cl R_R(\mathbf{E})$ is a hypersurface
and $\cl A _{R/(f)}(E)$ is Cohen--Macaulay of codimension $3$ (but not Gorenstein)  -- note the latter coincides with
the symmetric algebra $\cl S_{R/(f)}(E)$ since $\mathbf{E}$ is of linear type because it has
homological dimension $1$ with presentation $0\rar R\stackrel{\phi}{\lar} R^4\rar \mathbf{E}\rar 0$
having $I_1(\phi)$ of codimension $\geq 2$.
On the other hand, $\cl R _{R/(f)}(E)$ is almost Cohen--Macaulay, but not Cohen--Macaulay.

For the nodal cubic  $\cl R_R(\mathbf{E})$ is a codimension $2$ complete intersection and $\cl A _{R/(f)}(E)$ is almost Cohen--Macaulay, but not Cohen--Macaulay.

\subsection{The normal module of a Cohen--Macaulay ideal of codimension $2$}

Let $J\subset R$ be an ideal in a Noetherian ring $R$ with a given free presentation
$$R^m \stackrel{\phi}{\lar} R^n \lar J\rar 0.$$
This induces a free presentation of the conormal module $J/J^2$ over $R/J$ which gives
the normal module as a syzygy module
$$0\rar (J/J^2)^*\lar {(R/J)^n}^* \stackrel{\bar{\phi}^t}{\lar}  {(R/J)^m}^*.$$
We take $E:=(J/J^2)^*$ and  $J$-deform $E$ to its inverse image $\mathbf{E}$ in ${R^n}^*$
by the natural surjection ${R^n}^*\surjects {(R/J)^n}^*$.

Thus, in an explicit way, $\mathbf{E}$ is the sum of the $R$-submodule of ${R^n}^*=\mathbf{F}$ generated by the liftings of the generators of $E$ as a submodule of ${(R/J)^n}^*$ and $J \mathbf{F}$.

\begin{Remark}\rm We note that the case study of the previous section is in principle subsumed under the present case study. Quite generally, if $A$ is an algebra essentially of finite type over a field $k$, the module of derivations ${\rm Der}_k(A)$ is identified with the $A$-dual of the conormal module $\mathbb{D}_A/\mathbb{D}_A^2$
of the diagonal embedding  $J:=\mathbb{D}_A\subset R:=A\otimes _kA$.
\end{Remark}

In the case where $R$ is a local or standard graded ring (with $J$ homogeneous), in general collecting generators from both summands does not lead to a minimal set of generators of $\mathbf{E}$.
It might actually happen that the lifting of a minimal set of generators of the normal module already generates $\mathbf{E}$, as was the case of the module of logarithmic derivations.

Since the normal module of a complete intersection is free, its Aluffi algebra is of no interest.
The next simplest case is that of a codimension $2$ perfect ideal.
The next theorem is a pointer for the boundaries of a general theory.

One of its assertions hinges on the following notion.
Let $R:=k[\xx]=k[x_1,\ldots,x_d]$ stand for a polynomial ring over a field $k$ and let $M$ denote a linearly presented $R$-module generated by $m$ elements, with presentation matrix $\Phi$.
Then the presentation ideal $I_1([\yy]\cdot \Phi)$ of the symmetric algebra $\mathcal{S}_R(M)$ can be written in the form $I_1([\xx]\cdot B(\Phi))$, where $[\yy]=[y_1\cdots y_m]$ are presentation variables and $B(\Phi)$ denotes the transposed Jacobian matrix of the generators of $I_1([\yy]\cdot \Phi)$ with respect to the $\xx$-variables.
Note that the entries of this matrix are linear forms in the polynomial ring $k[\yy]$.

The matrix $B(\Phi)$ is often called the {\em Jacobian dual matrix} of $M$ (or of its linear presentation).
By Cramer, the $d$-sized minors of $B(\Phi)$ belong to the defining ideal of the Rees algebra of $M$.
Typically, $d$ is too large a degree in order that these minors be minimal generators.
We will say that $M$ (or its Rees algebra) is of the {\em expected} fiber type provided the presentation ideal of the Rees algebra $\mathcal{R}_R(M)$ is $(I_1([\yy]\cdot \Phi), I_d(B(\Phi)))$.

\begin{Theorem}\label{Main_Normal}  Fix two integers $d\geq 1, n\geq 3$. Let $R:=k[\xx]=k[x_1,\ldots,x_d]$ be  a standard homogeneous polynomial ring over a field and let $\phi$ denote an $n\times (n-1)$ matrix $\phi$ whose entries are linear forms in $R$ {\rm (}possibly zero{\rm )}. Assume that the ideal   $J\subset R$ generated by the $(n-1)$-minors of $\phi$ has codimension $2$.
Then  a $J$-deformation $\mathbf{E}\subset \mathbf{F}={R^n}^*$ of the normal module $(J/J^2)^*$ satisfies the following properties:
\begin{enumerate}
\item[{\rm (i)}] The $R$-dual  $\mathbf{E}^*$ is a free module.

\item[{\rm (ii)}] Both  $(J/J^2)^*$ and  $\mathbf{E}$ specialize from the generic case of $\phi;$ in particular, $\mathbf{E}$ is the image of the map defined by the specialization of the Jacobian matrix of $J$ generic.
\item[{\rm (iii)}]    $\mathbf{E}$ is a linearly presented  module of projective dimension $1$ and minimally generated by $n(n-1)$ elements of degree $n-2$.
\item[{\rm (iv)}]   If $\mathbf{E}$ satifies the condition $G_d$ then
the Rees algebra of $\mathbf{E}$ is Cohen--Macaulay and of the expected fiber type.
\item[{\rm (v)}] $\mathbf{E}$ is a module of linear type if and only if it satisfies the condition $G_d$ and the inequality $n(n-2)\leq d-1$ holds.
\item[{\rm (vi)}] If the dimension of the symmetric algebra $\mathcal{S}_{R/J}((J/J^2)^*)$ of the normal module equals $n(n-1)$ then the dimension of its Aluffi algebra
$ \cl A_{R/J}((J/J^2)^*)$  equals $\ell (\mathbf{E})$, where the latter denotes the analytic spread of $\mathbf{E}$.
\end{enumerate}
\end{Theorem}
\demo
(i) We argue that $\mathbf{E}$ is an {\em ideal module} in the sense of {\rm\cite[Section 5]{ram1}}. For this, since $J$ has grade $2$, it suffices to show that $J\subset{\rm Ann}(\mathbf{F}/\mathbf{E})$ (\cite[Proposition 5.1 (c)]{ram1}).
But this is obvious since $J\mathbf{F}\subset \mathbf{E}$ by the construction of a deformation.

\medskip

(ii)  The degree of the generators of $J$ is $n-1$, hence the syzygies of the transpose $\phi^t$ modulo $J$ have initial degree bounded below by $n-2$. That the generating syzygies have exactly this degree is to be seen.

Consider preliminarily the fully generic case, i.e., where the entries of $\phi$ are indeterminate entries.
It is well-known that the ideal generated by the maximal minors of a non-square generic matrix is rigid \cite[Theorem 5.2]{Svanes}.
This means that the normal module coincides with the Jacobian module of the generators of $J$ (i.e., the submodule which is the image of the map defined by the Jacobian matrix of these generators considered modulo $J$).
Clearly then the normal module is minimally generated by $n(n-1)$ elements.

\smallskip

We next specialize to the given $\phi$ by mapping the entries of the generic matrix to the respective entries of $\phi$.
Namely, let $S:=R[X_{i,j}\,|\, 1\leq i\leq j\leq n-1]$ be the polynomial ring on the entries of the generic matrix extending further the coefficients to $R$.
Map $S\rar R$ via $X_{i,j}\mapsto \ell_{i,j}$, where $\phi=(\ell_{i,j})_{1\leq i\leq j\leq n-1}$.
Since the entries $\ell_{i,j}$  are linear, we may harmlessly assume that they generate the maximal ideal $(x_1,\ldots,x_d)$.
Therefore, the map $S\rar R$ is surjective, with kernel the regular $S$-sequence $\{x_{i,j}-\ell_{i,j}\,|\, 1\leq i\leq j\leq n-1\}$.
The latter is also a regular sequence modulo the ideal of maximal minors of the generic matrix.

Next one uses a property that holds for an arbitrary perfect codimension $2$ ideal $I\subset A$,  namely, such an ideal is a licci ideal. Since a licci ideal is strongly non-obstructed,  $I/I^2\otimes_A \omega_{A/I}$ is a maximal Cohen--Macaulay module, where $\omega_{A/I}$ is the canonical module of $A/I$.
Therefore, to show that the normal module specializes, it suffices to show this is the case by going modulo a nonzerodivisor $a\in A$ on both $A$ and $A/I$.
Now, $a$ is then also a regular element on $\omega:=\omega_{A/I}$ as the latter is a maximal Cohen--Macaulay module.
On the other hand, since $I/I^2\otimes_A \omega$ is a maximal Cohen--Macaulay module one has Ext$^1(I/I^2\otimes_A \omega,\omega)=0$.
Therefore, the exact sequence $0 \to \omega \to \omega \to \omega/a \omega \to 0$ induces an exact sequence
$$0 \to {\rm Hom}(I/I^2\otimes_A \omega, \omega) \stackrel{.a}{\to} {\rm Hom}(I/I^2\otimes_A \omega, \omega) \to {\rm Hom}\left(I/I^2\otimes_A \omega, \frac{\omega}{a\omega}\right) \to 0.$$
On the other hand, one has
\begin{eqnarray}\nonumber
{\rm Hom}(I/I^2\otimes_A \omega,\omega)&\simeq & {\rm Hom}(I/I^2,A)\\ \nonumber
{\rm Hom}(I/I^2\otimes_A \omega, \omega/a\omega)&\simeq & {\rm Hom}((I,a)/(I^2,a),A/(a)).
\end{eqnarray}

Therefore, one gets

\begin{eqnarray}\nonumber
((I,a)/(I^2,a))^*={\rm Hom}((I,a)/(I^2,a),A/(a)) &\simeq & {\rm Hom}(I/I^2,A)/a\,{\rm Hom}(I/I^2,A)\\ \nonumber
&=& (I/I^2)^*\otimes_A A/(a),
\end{eqnarray}
i.e., the normal module specializes.
(We are indebted to Bernd Ulrich for communicating to us this result and its proof).

Now let $\boldsymbol\Theta$ denote the Jacobian matrix over $S$ of the maximal minors in the fully generic case of $\phi$.
Thus, the normal module in the generic case is the submodule generated by the column vectors of $\boldsymbol\Theta$ considered modulo the generic maximal minors.
Letting $\widetilde{\boldsymbol\Theta}$ stand for the matrix over $R$ specialized from $\boldsymbol\Theta$, we have shown that the normal module $(J/J^2)^*$ over $R/J$ is  the submodule image of $\widetilde{\boldsymbol\Theta}$ read modulo $J\subset R$.
As such, in particular it is generated by  $n(n-1)$ elements.

By definition, to get a deformation $\mathbf{E}$ of the normal module $\mathbf{F}=R^n$ one takes the $R$-submodule generated by the image of  $\widetilde{\boldsymbol\Theta}$ and $J\mathbf{F}$. Note however that the product of $^t\phi$ by any column vector of $\widetilde{\boldsymbol\Theta}$ not only belongs to $J$, it is actually a $k$-linear combination of its set of generators.
This shows that $J\mathbf{F}$ is already contained in the $R$-submodule generated by the image of  $\widetilde{\boldsymbol\Theta}$, hence
 $\mathbf{E}$ is generated by the latter.

\medskip

(iii) Consider again the generic case, as in (ii).
Set $S:=k[X_{i,j}\,|\, 1\leq i<j\leq n-1]$.
Let $\boldsymbol\Theta$ denote the Jacobian matrix of the maximal minors of the matrix ${\phi}:= (X_{i,j})_{1\leq i<j\leq n-1}$.
We will show that the $S$-module $\mathcal{E}$ image of the corresponding Jacobian map has projective dimension $1$ and then apply generic perfection to get an explicit resolution of the $R$-module obtained as the image of the Jacobian map by specializing the entries of ${\phi}$ to the respective entries of the given matrix in the statement of the theorem.

It will be convenient to rewrite $\boldsymbol\Theta$ by changing basis of the free ambient $S^n$ and generators of $\mathcal{E}$.
Thus doing, using the same notation, one can write $\boldsymbol\Theta$ as the concatenation of $n-1$ matrices
$$\boldsymbol\Theta=\left(\;\Theta_1\,|\ldots|\,\Theta_{n-1}\;\right)$$
where

$$\Theta_j=\left(\begin{array}{cccc}
0&\delta^{(2)}_{1,j}&\ldots&\delta^{(n)}_{1,j}\\
\delta^{(1)}_{2,j}&0&\ldots&\delta^{(n)}_{2,j}\\
\vdots&\vdots&\ddots&\vdots\\
\delta^{(1)}_{n,j}&\delta^{(2)}_{n,j}&\ldots&0
\end{array}
\right).
$$
Here for each $1\leq i\leq n$, the symbol $\delta^{(i)}_{k,j}$ denotes the signed cofactor of the $(k,j)$-entry in the $(n-1)\times (n-1)$ matrix $\varphi_{(i)}$ obtained form $\phi$ by omitting its $i$th row.

We now describe a set of natural syzygies of $\Theta$.
Namely, let $C_j$ denote the $j$th column of $\phi$ and let $\phi^{(j)}$ stand for the $n\times (n-2)$-submatrix of $\phi$ by omitting $C_j$.
The following relations are immediately seen by Cramer's rule:

\begin{equation}\label{identidades}
\Theta_j\cdot C_j=(D_1\ldots,D_n)^t\quad\quad \mbox{e}\quad\quad \Theta_j\cdot \varphi^{(j)}=\boldsymbol0
\end{equation}
where $D_i$ stands for the signed $n-1$-minor of $\phi$ omitting the $i$th row.

The matrix of syzigies corresponding to these relations can be written as the concatenation of two matrices, the first an
$n(n-1)\times (n-1)(n-2)$ block diagonal matrix and the second an
$n(n-1)\times (n-2)$ matrix in semi-block diagonal form:

\begin{equation}\label{true_syzygies}
\Phi=
\left(\begin{array}{cccccc|cccc}
\phi^{(1)}&&&&& &-C_1 &&&\\
&\phi^{(2)}&&&& &\kern5pt C_2 &-C_2 &&\\
&&\phi^{(3)}&&&& &\kern5pt C_3 &&\\ %[-15pt]
&&&\ddots && &&&\ddots &\\
&&&&\phi^{(n-2)}&&&&& -C_{n-2}\\
&&&&&\phi^{(n-1)}&&&& C_{n-1}
\end{array}
\right),
\end{equation}
where the vacant spaces are filled with zero entries.

Therefore, one has the following complex of free modules
$$0\rar S^{n(n-2)}\stackrel{\Phi}\lar S^{n(n-1)}\stackrel{\boldsymbol\Theta}\lar S^{n}.$$

To prove exactness we draw upon the well-known Buchsbaum--Eisenbud acyclicity criterion. Thus, first compute the ranks of $\Theta$ and of $\Phi$.
As for $\boldsymbol\Theta$ it is well-known that the maximal minors of $\phi$ are algebraically independent over $k$.
Therefore, the rank of their Jacobian matrix, and hence of the present matrix $\boldsymbol\Theta$, is $n$.

\smallskip

In the case of $\Phi$, we note that from (\ref{true_syzygies}), it admits the following $n(n-2)\times n(n-2)$ submatrix in semi-diagonal block shape:

$$\Psi=\left(\begin{array}{cc|cc|cc|cc|cccc}
\phi^{(1)}_{i}&C_{1,i}&&&&&&&\\
\hline
&C_{2,i}&\phi^{(2)}_{i}&C_{2,i}&&&&&\\
\hline
&&&&\ddots&&&&&\\
\hline
&&&&&C_{n-2,i}&\phi^{(n-2)}_i&C_{n-2,i}\\
\hline
&&&&&&&\ast&\psi
\end{array}\right).
$$
Here $\ast$ denotes a certain $(n-2)\times 1$ submatrix of $C_{n-1}$
and, for $1\leq i\leq n$ and $1\leq j\leq n-2,$
$\phi^{(j)}_i$ and  $C_{j,i}$ stand, respectively, for the submatrices obtained from $\phi^{(j)}$ and $C_{j}$  omitting their respective $i$th rows -- note that thus $(\phi^{(j)}_{i}|C_{j,i})$ is a square matrix of order $n-1$.
In addition,  $\psi$ is an arbitrary $(n-2)\times (n-2)$ submatrix of $\phi^{(n-1)}$ with non-vanishing determinant -- such matrices exist since  $J=I_{n-1}(\phi)\subset I_{n-2}(\phi^{(n-1)})$.

Then $$\det(\Psi)=\det(\phi^{(1)}_{i}|C_{1,i})\cdots\det(\phi^{(n-2)}_i|C_{n-2,i})\det(\psi)=\pm (D_{i})^{n-2}\det(\psi).$$
where $D_i$ denotes as earlier the $(n-1)$-minor of $\phi$ striking out the $i$th row.

\smallskip

It remains to check the codimensions of the relevant Fitting ideals.
Since $\Theta $ has rank $n$, clearly $I_n(\Theta)$ has codimension at least $1$.
As for $I_{n(n-2)}(\Phi)$, since $\psi$ is an arbitrary $(n-2)\times (n-2)$ submatrix of $\phi^{(n-1)}$ with non-vanishing determinant, we observe that the above determinant actually encapsulates the following ideal inclusion
\begin{equation}\label{inclusion}
((D_1)^{n-2},\ldots,(D_n)^{n-2})\cdot I_{n-2}(\phi^{(n-1)})\subset I_{n(n-2)}(\Phi).
\end{equation}

Since $I_{n(n-2)}(\Phi)$ contains the product of two ideals each of codimension at least $2$ then ${\rm cod}(I_{n(n-2)}(\Phi))\geq 2$, as required.

\smallskip

To conclude, we observe that the cokernel of $\boldsymbol\Theta$ is a perfect $S$-module. Indeed, its annihilator has same codimension as the Fitting ideal $I_n(\boldsymbol\Theta)$.
But in the generic case, the maximal minors of the Jacobian matrix have same codimension as the ideal $J:=(D_1,\ldots,D_n)\subset S$ of $n$-minors of $\phi$.
To see this, one has on one hand that $\boldsymbol\Theta$ has rank $2$ modulo the ideal $J$ of maximal minors, which says that $I_n(\boldsymbol\Theta)\subset J$ (recall that $n\geq 3$).
On the other hand, in the opposite direction, one has $((D_1)^{n-2},\ldots,(D_n)^{n-2})\subset I_n(\boldsymbol\Theta)$, hence the two ideals have the same codimension (in fact, share the same radical).
This latter inclusion clearly holds after specialization, hence $I_n(\widetilde{\boldsymbol\Theta})$ has codimension at least $2$, where $\widetilde{\boldsymbol\Theta}$ denotes the specialized Jacobian matrix.
By generic perfection (see \cite[Theorem 3.5]{BHbook}) we get a free resolution

 $$0\rar R^{n(n-2)}\stackrel{\widetilde{\Phi}}\lar R^{n(n-1)}\stackrel{\widetilde{\boldsymbol\Theta}}\lar R^{n},$$
 where $\widetilde{\Phi}$ denotes the  specialization of $\Phi$.

 \smallskip

To conclude, note that the image of $\widetilde{\boldsymbol\Theta}$ is the deformation $\mathbf{E}$ by the previous item.
Therefore, this proves that  $\mathbf{E}$ is a linearly presented $R$-module of projective dimension $1$. Moreover, since all entries of the matrices in the above free resolution are homogeneous of positive degree, the assertion about the minimal number of generators of $\mathbf{E}$ also follows.

\medskip

(iv)  Using (iii), this reads off \cite[Proposition 4.11]{ram1}.

\medskip

(v) We first translate the condition $G_d$ in terms of Fitting ideals in the present context.
For convenience, write $\Phi$ (instead of $\widetilde{\Phi}$) for the linear presentation matrix of $\mathbf{E}$ found in item (iii).

Applying the equivalence explained in \cite[Section 3, p. 615]{ram1} to our context one sees that $\mathbf{E}$ satisfies $G_d$ if and only if $I_t(\Phi)\geq n(n-2)-t+2$ for $n(n-2)-(d-2)\leq t\leq n(n-2)$.
But since we are assuming the inequality $n(n-2)\leq d-1$ then the leftmost inequality of the above interval for $t$ implies that the interval becomes $1\leq t\leq n(n-2)$.
In other words, $\mathbf{E}$ actually satisfies $G_{\infty}$ (also called $F_1$).
But a module of projective dimension one satisfying   $G_{\infty}$ is of linear type (\cite[Theorem 3.4]{SV1}).

Conversely, if $\mathbf{E}$ is of linear type it satisfies $G_d$ for even more reason and besides, the inequality $n(n-2)\leq d-1$ is equivalent to the following one:
$$n(n-1)=\mu(\mathbf{E})\leq \dim R+rk(E)-1= d+n-1,
$$
the latter being satisfied by the  $G_{\infty}$ condition on the nose, as read in terms of number of generators of $\mathbf{E}$.

\medskip

(vi) Recall that
\begin{eqnarray}\nonumber
\cl A_{R/J}((J/J^2)^*)&\simeq & \mathcal{R}_R(\mathbf{E})/(J,J\mathbf{F})\\ \nonumber
&\simeq & R[\yy]/(J,I_1(\yy\cdot \Phi), \mathfrak{A}_{\geq 2}, I_1(\yy\cdot\Psi)\\ \nonumber
&\simeq & \frac{\mathcal{S}_{R/J}((J/J^2)^*)}{\widetilde{\mathfrak{A}_{\geq 2}}}
\end{eqnarray}
where $\yy=\{y_{i,j} \,|\, 1\leq i\leq n,1\leq j\leq n-1\}$,  $\mathfrak{A}_{\geq 2}$ denotes the Rees equations of $\mathbf{E}$ of $\yy$-degree at least $2$, whereas the tilde over it has the obvious meaning, and finally $\Psi$ denotes the content matrix of the inclusion $J\mathbf{F}\subset \mathbf{E}$ (a matrix with entries on $R$).

By assumption, $\dim \mathcal{S}_{R/J}((J/J^2)^*)=n(n-1)$.

Reading this value backwards in the above string of isomorphisms, the defining ideal $(J,I_1(\yy\cdot \Phi), I_1(\yy\cdot\Psi)$ of $\mathcal{S}_{R/J}((J/J^2)^*)$ has codimension $d+n(n-1)-n(n-1)=d$ on the ambient polynomial ring $R[\yy]$.
Since it is wholly contained in the extension $(\xx)R[\yy]$ of the maximal ideal $(\xx)$ of $R$, it must be primary to $(\xx)R[\yy]$.
If now $P\subset R[\yy]$ is any prime ideal containing the entire Aluffi defining ideal $(J,I_1(\yy\cdot \Phi), I_1(\yy\cdot\Psi), \mathfrak{A}_{\geq 2})$ it must contain the ideal $(\xx,\mathfrak{I})$, with $\mathfrak{I}$ standing for the ideal generated by the bigraded piece of $\mathfrak{A}_{\geq 2}$ of bidegrees $(0,\geq 2)$.
But $\mathfrak{I}$ is the defining ideal of the $k$-subalgebra $k[\mathbf{E}]\subset \mathcal{R}_R(\mathbf{E})$ over $k[\yy]$ -- i.e., the special fiber of $\mathcal{R}_R(\mathbf{E})${\rm )}.
In particular, $(\xx,\mathfrak{I})$ is a prime ideal and must be the radical of the Aluffi defining ideal.

It follows that $\dim \cl A_{R/J}((J/J^2)^*)=\dim R[\yy]/(\xx,\mathfrak{I})=\dim k[\yy]/\mathfrak{I}$, as was to be shown as the analytic spread is the dimension of the special fiber.
\qed

\medskip

We conclude the section with a couple of remarks and a few easy examples in low dimensions (computed with the help of \cite{Macaulay})

\begin{Remark}\rm
The dimension theoretic hypothesis of item (vi) above is quite natural and is verified in several situations.
Since $n(n-1)$ is the number of minimal generators of the normal module and $n(n-1)\geq d=\dim R$ is the standing proviso, then this hypothesis translates into saying (at least in the case that $J$ is generically a complete intersection) that the normal module is of {\em Valla type} in the sense of an analogue for modules of the notion introduced for ideals in \cite[Section 3]{Salvador} (see also \cite[Section 2]{s1}).
On the other hand, one could easily show that this hypothesis holds in general provided it would hold in the generic case. In the latter situation, the assumption also means that the symmetric algebra of the normal module has the expected dimension, i.e., same as that of the Rees algebra of the module.
Unfortunately, some experimentation shows that the normal module in the generic case may fail to satisfy the required condition $(F_0)$, for $n\geq 6$.
Since in the generic case the ideal $J$ is rigid, the first syzygies of the normal module turn out to be the $k$-derivations of the algebra $R/J$.
Thus, the question becomes one of estimating the codimensions of the ideals of minors of the map whose image is the module of $k$-derivations of $R/J$ -- a rather fundamental corner of rigidity.
\end{Remark}

\begin{Remark}\rm
A natural question arises as to whether one can treat similarly the case of the ideal $J$ of maximal minors of an $n\times m (m < n)$ matrix having codimension $n-m+1$, given that the latter is also rigid. A technical difficulty arises at the outset due to the fact that a deformation $\mathbf{E}\subset \mathbf{F}$ of the corresponding normal module will require quite a bit of minimal generators coming from the submodule $J\mathbf{F}$, hence the numbers will be harder to control.
\end{Remark}

\begin{Example}\rm (The coordinate axes of $\pp^2$)

Let $J=(xy,xz,yz)\subset R=k[x,y,z]$.
Here $\mathbf{E}$ satisfies the condition $G_3$, hence by (iv) the presentation ideal of its Rees algebra has the expected fiber type, with $k[\mathbf{E}]$ defined by the  determinant of the Jacobian dual matrix. Clearly, then $\mathbf{E}$ is not of linear type -- and indeed, $n(n-2)=3=d$ which is in accordance with (v).
A computation yields that the Aluffi algebra $\mathcal{A}$ has indeed dimension $5$, in accordance with (vi).
Neither $\mathcal{A}$ nor the Rees algebra of the normal module is a Cohen--Macaulay algebra.
\end{Example}

\begin{Example}\rm (Four cubics defining the Magnus involution in $\pp^3$)

Let $J=(xyz,xyw, xzw,yzw)\subset R=k[x,y,z,w]$.
Here $\mathbf{E}$ does not satisfy the condition $G_4$; in particular, by (v) it is not of linear type. The minimal defining equations of $k[\mathbf{E}]$ are quadrics and cubics, while the rank of the Jacobian dual matrix $B(\Phi)$ is $d=4$.
The Rees algebra of $\mathbf{E}$ is still of  fiber type, but $I_d(B(\Phi))$ is no longer minimal generators -- although, as always, $I_d(B(\Phi))$ is contained in the presentation ideal of the special fiber algebra $k[\mathbf{E}]$.
A computation corroborates that the unmixed part of the Aluffi algebra has the form stated in the proof of (vi); in particular, $\dim \mathcal{A}=\ell(\mathbf{E})=7$.
Neither $\mathcal{A}$ nor the Rees algebra of the normal module is a Cohen--Macaulay algebra.
\end{Example}

\begin{Example}\rm (The defining equations of the twisted cubic)

Let $J=(xz-y^2,xw-yz,yw-z^2)\subset R=k[x,y,z,w]$.
Here $\mathbf{E}$  satisfies the condition $G_4$ and $n(n-2)=3=d-1$. Therefore, $\mathbf{E}$ is of linear type.
By the proof of (vi) in this situation, $\ell(\mathbf{E})$ is maximal and  $\dim \mathcal{A}=\ell(\mathbf{E})=6$.
$\mathcal{A}$ is not Cohen--Macaulay, but the Rees algebra of the normal module is a Cohen--Macaulay domain.
\end{Example}

%{\sc(\% $J$ should theoretically function for $\mathbf{E}$ as a sort of replication of a Bourbaki ideal)}

%\begin{Remark}\rm
%A Bourbaki ideal of $\mathbf{E}$ can be obtained as follows: let $\mathbf{E}'\subset \mathbf{E}$ be a submodule generated by $\rk (\mathbf{E})-1= \mu(J)-1$ general columns. Then the $R$-module $N:=\mathbf{E}/\mathbf{E}'$ has a free resolution
%$$0\rar R^{m-1}\stackrel{\Phi}{\lar} R^m \lar N\rar 0,$$
%where $m=\mu(\mathbf{E})$.
%Then $N$ identified with the ideal generated by the maximal minors of $\Phi$ is a Bourbaki ideal of $\mathbf{E}$.

%(Note that $\Phi$ has the identity matrix of size $\mu(J)-1$ as a diagonal block (up to elementary transformations)).
%\end{Remark}

\end{document}